\documentclass{article}

\input{header}
\usepackage[a4paper,top=3cm,bottom=2cm,left=3cm,right=3cm,marginparwidth=1.75cm]{geometry}

\makeatletter
\newcommand\footnoteref[1]{\protected@xdef\@thefnmark{\ref{#1}}\@footnotemark}
\makeatother

\usepackage{amsmath,amsfonts,amssymb,amsthm,dsfont,natbib,algorithm,algorithm}
\usepackage{authblk}

\begin{document}
\title{Stochastic Scale Invariant Power Iteration for KL-divergence Nonnegative Matrix Factorization}
\author[1]{Cheolmin Kim}
\author[2]{Youngseok Kim}
\author[1]{Diego Klabjan}
\affil[1]{Department of Industrial Engineering and Management Sciences, Northwestern University}
\affil[2]{Department of Statistics, University of Chicago}
\date{\vspace{-5ex}}

\maketitle

\begin{abstract}
We introduce a mini-batch stochastic variance-reduced algorithm to solve finite-sum scale invariant problems which cover several examples in machine learning and statistics such as principal component analysis (PCA) and estimation of mixture proportions. The algorithm is a stochastic generalization of scale invariant power iteration, specializing to power iteration when full-batch is used for the PCA problem. In convergence analysis, we show the expectation of the optimality gap decreases at a linear rate under some conditions on the step size, epoch length, batch size and initial iterate. Numerical experiments on the non-negative factorization problem with the Kullback-Leibler divergence using real and synthetic datasets demonstrate that the proposed stochastic approach not only converges faster than state-of-the-art deterministic algorithms but also produces excellent quality robust solutions.
\end{abstract}

% main body
%---------------------------------------------%
\section{Introduction}
%---------------------------------------------%
We study a class of optimization problems called \textit{finite-sum scale invariant problems} of the form
\begin{equation} \label{prob:main}
\begin{aligned} 
    & \max & \quad & f(x) =  \textstyle \frac{1}{n} \sum_{i=1}^n f_i(x) \quad \\
    & \textrm{subject to} & \quad & x \in \partial \mathcal{B}_d \triangleq \{ x \in \mathbb{R}^d \, | \, \| x \| = 1 \} 
\end{aligned}
\end{equation}
where $f_i$ are scale invariant functions of the same type, i.e. $f_i$ are either multiplicatively scale invariant such that $f_i(cx) = u(c) f_i(x)$ for the same multiplicative factor $u(c)=|c|^p$ or additively scale invariant satisfying $f_i(cx) = f_i(x)+v(c)$ with the same additive factor $v(c) = \log_a |c|$. All norms are 2-norms unless indicated otherwise. 
The scale invariant problem covers interesting problems in machine learning and statistics such as $L_p$-norm kernel PCA \citep{kim2019simple} and estimation of mixture proportions \citep{kim2018fast}, to name a few. Moreover, as studied in \cite{kim2019scale}, more examples such as independent component analysis (ICA) \citep{hyvarinen2004independent,hyvarinen1999fast}, Gaussian mixture models (GMM), Kullback–Leibler divergence non-negative matrix factorization (KL-NMF) \citep{fevotte2011algorithms, lee2001algorithms, wang2013nonnegative} and the Burer-Monteiro factorization of semi-definite programs \citep{erdogdu2018convergence} can be formulated to extended settings of \eqref{prob:main}.

If the objective function $f$ is twice differentiable on an open set containing $\partial \mathcal{B}_d$, the scale invariant problem \eqref{prob:main} can be locally viewed as the leading eigenvector problem in the sense that a stationary point $x^*$ is an eigenvector of $\nabla^2 f(x^*)$. Moreover, if the Lagrange multiplier $\lambda^*$ satisfying $\lambda^* x^* = \nabla f(x^*)$ is greater than $\bar{\lambda}$, the largest absolute eigenvalue of $\nabla^2 f(x^*)(I-x^*(x^*)^T)$, the stationary point $x^*$ is a local maximum.
Due to this eigenvector property, the scale invariant problem can be efficiently solved by a general form of power iteration called \textit{scale invariant power iteration (SCI-PI)} \citep{kim2019scale} which repeats
\begin{equation} \label{eq:sci-pi}
    x_{k+1} \leftarrow \frac{\nabla f(x_k)}{\| \nabla f(x_k) \|}.
\end{equation}

Interestingly, the convergence of SCI-PI generalizes that of power iteration. If $x_0$ is initialized close to a local optimum $x^*$, the optimality gap $1-(x_k^Tx^*)^2$ linearly converges to zero at an asymptotic rate of $(\bar{\lambda}/\lambda^*)^2$. 
For PCA \citep{jolliffe2002principal}, this rate specializes to $(\lambda_1/\lambda_2)^2$ where $\lambda_1$ and $\lambda_2$ are the first and the second eigenvalues of the covariance matrix $\frac{1}{n} \sum_{i=1}^n a_i a_i^T$ constructed by data vectors $a_i$. This convergence analysis shows that SCI-PI not only has a general form of power iteration but also extends the attractive local linear convergence property of power iteration.

The convergence analysis of power iteration for PCA is analogous to that of gradient descent for strongly convex optimization problems. Due to this analogy, many variants have been developed for power iteration such as noisy \citep{hardt2014noisy}, momentum \citep{xu2018accelerated}, coordinate-wise \citep{lei2016coordinate}, online \citep{garber2015online, warmuth2008randomized, boutsidis2015online, nie2016online} and stochastic \citep{oja1982simplified, oja1985stochastic, shamir2016convergence, arora2012stochastic, arora2013stochastic} power methods.
In particular, using the stochastic variance-reduced gradient technique \citep{johnson2013accelerating}, stochastic variance-reduced PCA algorithms \citep{shamir2015stochastic, xu2018accelerated,kim2019stochastic} have reduced the total runtime to obtain an $\epsilon$-optimal solution from $\mathcal{O} ( dn (\frac{\lambda_1}{\lambda_1-\lambda_2}) \, \text{log} \frac{1}{\epsilon} )$ to ${\mathcal{O}} ( d ( n+ (\frac{\lambda_1}{\lambda_1-\lambda_2} )^2 \big) \, \text{log} \frac{1}{\epsilon} )$. This decoupling of the sample size $n$ from the eigen-gap $1-\lambda_2/\lambda_1$ is advantageous in a large scale setting where $n$ is relatively larger $\lambda_1/(\lambda_1-\lambda_2)$.

In this work, we introduce a stochastic algorithm with mini-batch variance-reduced gradients called \textit{Stochastic Scale Invariant Power Iteration (S-SCI-PI)} to solve finite-sum scale invariant problems \eqref{prob:main} and provide a convergence analysis for it. 
While stochastic algorithms have been extensively studied in statistics and machine learning, to our best knowledge, there has been no work in the literature that develops a stochastic algorithm for finite-sum scale invariant problems. With stochastic variance-reduced gradients, we can improve the total runtime by decoupling the sample size $n$ from the eigen-gap $1 - \bar{\lambda} / \lambda^*$, so we can efficiently solve a complex matrix problems like KL-NMF when $n \gg \lambda^* / (\lambda^* - \bar{\lambda})$. Although some stochastic variance-reduced algorithms have been developed for nonconvex matrix problems such as nonconvex low-rank matrix recovery \citep{wang2017unified} and Frobenius-norm NMF \citep{kasai2018stochastic}, no stochastic variance-reduced algorithm with convergence guarantee has been developed for the KL-NMF problem. \citep{serizel2016mini} introduced some stochastic mini-batch updates for the KL-NMF problem but it lacks a convergence analysis.

S-SCI-PI is a generalization of SCI-PI and VR Power \citep{kim2019stochastic}. Using the observation that finite-sum scale invariant problems can be locally seen as the PCA problem, we adapt VR Power, which has a provable optimal runtime for any batch size. Specifically, we adjust the scaling factor for the full-gradient depending on the degree of scale invariance of the objective function. In convergence analysis, we prove linear convergence of the expected optimality gap, which is expressed as the ratio of two expectation terms.
As in the analysis of VR Power, the convergence rate of S-SCI-PI depends on the eigenvalues of the Hessian at the solution. However, the analysis of VR Power is not easily transferred since S-SCI-PI has an additional error term arising from the difference of the Hessians between the iterate and the optimal solution. We provide a condition on the step size, epoch length, batch size and initial iterate, which ensures that this error is not increasing in the course of the algorithm and that the expected optimality gap converges at a linear rate.

%In brief, our stochastic S-SCI-PI algorithm can be understood as the locally VR-Power-like algorithm, by decomposing the objective function into the sum of the local quadratic term and the additional error term. To bound this additional error term, we introduce mild conditions on step size, batch size and the initial iterate. This ensures that the error is not increasing in the course of the algorithm and guarantees the linear convergence of the expected optimality gap.

%The first error term is attributed to the difference of the Hessians between the iterate and the optimal solution. To control this error, we derive a condition on the step size, batch size and initial iterate, which ensures that the error is not increasing in the course of the algorithm. On the other hand, the second error occurs from the stochastic sampling of gradient. 
%Using recursion, we develop a compact representation of the optimality gap. Moreover, we derive a condition on the step size, epoch length, batch size and initial iterate, which ensures the linear convergence of the expected optimality gap.

In numerical experiments, we report a study on the KL-NMF problem. As reported in \citep{kim2019scale}, the KL-NMF problem can be decomposed into scale invariant subproblems where each subproblem considers a column of a matrix. By alternatively applying S-SCI-PI to these subproblems, we compute an optimal solution to the KL-NMF problem. Experiments on synthetic and real datasets demonstrate that the proposed stochastic approach not only converges faster than state-of-the-art deterministic algorithms but also produces robust solutions under random initialization.

Our work has the following contributions.
\iffalse
$\bullet $ We propose a stochastic algorithm (S-SCI-PI) for solving the finite-sum scale invariant problem. The algorithm adapts the stochastic variance-reduced gradient technique and adjusts the scaling factor of full-gradients depending on the order of scale invariance.

$\bullet $  We provide a convergence analysis for S-SCI-PI. Deriving compact representations of error terms, we prove linear convergence of S-SCI-PI such that the expected optimality gap decreases at a linear rate under some conditions on the initial iterate, epoch length, batch size and an additional condition.

$\bullet $ We provide experiments showing that SCI-PI converges faster than the state-of-the-art deterministic algorithms for KL-NMF and its subproblem.
\fi
\begin{itemize}
    \item We propose the stochastic algorithm S-SCI-PI to solve finite-sum scale invariant problems. The algorithm adapts the stochastic variance-reduced gradient technique by adjusting the scaling factor of full-gradients depending on the order of scale invariance.
    \item We provide a convergence analysis for S-SCI-PI. Deriving compact representations of error terms, we prove linear convergence of S-SCI-PI where the expected optimality gap decreases at a linear rate under some conditions on the step size, epoch length, batch size and initial iterate.
    \item We introduce new stochastic approach to solve the KL-NMF problem. Computational experiments show that our approach converges faster than state-of-the-art deterministic KL-NMF algorithms.
\end{itemize}
The paper is organized as follows. We present the algorithm in Section~\ref{sec:algorithm} and provide the convergence analysis in Section~\ref{sec:convergence-analysis}. We introduce the KL-NMF problem and its decomposition to scale invariant problems in Section~\ref{sec:KL-NMF}. We discuss some implementation issues in Section~\ref{sec:practical-considerations}. The experimental results on real and synthetic datasets are followed in Section~\ref{sec:experiment}.
%---------------------------------------------%
\section{Algorithm}
\label{sec:algorithm}
%---------------------------------------------%
Before presenting the algorithm, we first introduce some notations.
For the scale invariant objective function $f$ in \eqref{prob:main}, we let $p$ be the degree of scale invariance. If $f$ is multiplicatively scale invariant, $p$ is the order of the multiplicative factor $u(c)=|c|^p$. On the other hand, for additively scale invariant function, let $p=0$.
We denote the $k$-th coordinate of the gradient $\nabla f$ as $\nabla_k f(x)$. 
For a mini-batch sample $S \subset [n] \triangleq \{1,2,\cdots,n\}$, we define a stochastic function $f_{S} = \sum_{l \in S} {f_l}/{|S|}$. 
%We present the scale invariant multiplicative case and point out the changes needed for the additive case later. 

Following the stochastic variance-reduced gradient technique \citep{johnson2013accelerating}, our algorithm has a two-loop structure. 
At the start of each inner-loop, we compute the full gradient $\tilde{g}_T$ at the outer iterate $\tilde{x}_T$ and use this gradient information to construct a stochastic variance-reduced gradient $g_t$ at the inner iterate $x_t$. In order to derive a stochastic variance-reduced gradient at $x_t$ using the full gradient at $\tilde{x}_T$, we decompose $x_t$ as
\begin{equation*}
    x_t = \frac{x_t^T\tilde{x}_T}{\| \tilde{x}_T \|^2} \tilde{x}_T + x_t - \frac{x_t^T\tilde{x}_T}{\| \tilde{x}_T \|^2} \tilde{x}_T.
\end{equation*}
In the above equation, the first component is the projection of $x_t$ onto $\tilde{x}_T$ while the second part represents the orthogonal component of $x_t$ with respect to $\tilde{x}_T$. Since $\nabla f$ satisfies $c \nabla f(cx) = |c|^p \nabla f(x)$ \cite[Proposition 3]{kim2019scale}, assuming $x_t^T \tilde{x}_T > 0$, we can compute the exact gradient at the first component as
\begin{equation}
    \label{deterministic-gradient}
    \nabla f \left( \frac{x_t^T\tilde{x}_T}{\| \tilde{x}_T \|^2} \tilde{x}_T \right) = \frac{|x_t^T\tilde{x}_T|^{p-1}}{\| \tilde{x}_T \|^{2(p-1)}} \nabla f(\tilde{x}_T) = \alpha_t \tilde{g}_T
\end{equation}
where
\begin{equation*}
    \alpha_t = \frac{|x_t^T\tilde{x}_T|^{p-1}}{\|\tilde{x}_T\|^{2(p-1)}}.
\end{equation*}
To approximate the difference of gradients at $x_t$ and $(x_t^T\tilde{x}_T)\tilde{x}_T/\| \tilde{x}_T \|^2$, we use a stochastic sample $S_t \subset [n]$ of size $s$, which results in a stochastic variance-reduced gradient $g_t$ at $x_t$ as
\begin{equation*}
    \textstyle g_t = \alpha_t \tilde{g}_T + \frac{1}{s} \sum_{l \in S_t} \left[ \nabla f_l(x_t) - {\alpha_t} \nabla f_l(x_0) \right].
\end{equation*}
Since $g_t$ is subject to stochastic error, we introduce a step size $\eta \in (0,1]$ to control the progress of the algorithm depending on the variance of $g_t$. Using the step size $\eta$, we derive the following update rule
\begin{equation*}
    \textstyle x_{t+1} \leftarrow (1-\eta) x_{t} + \eta {g_t}/{\| x_t \|^{p-2}}.
\end{equation*}
Note that we divide $g_t$ by $\| x_t \|^{p-2}$ to match its scale with $x_t$ since $\nabla^2 f(cx) x= (p-1) c^{p-2} \nabla f(x)$ for $c > 0$ \cite[Proposition 3]{kim2019scale}.

Summarizing all the above, we obtain Algorithm~\ref{alg:Stochastic-SCI-PI}.
%---------------------------------------------%
\begin{algorithm}[ht]
   \caption{Stochastic SCI-PI (S-SCI-PI)}
   \label{alg:Stochastic-SCI-PI}
   \normalsize
\begin{algorithmic}
   \STATE {\bfseries Parameter:} step size $\eta \in (0,1]$, batch size $s$, epoch length $m$ \\
   \STATE randomly initialize outer iterate $\tilde{x}_0 \in \partial \mathcal{B}_d$ \\
   \FOR{$T=0,1,\ldots$}
      \STATE $x_0 \leftarrow \tilde{x}_{T}$, $\tilde{g}_T \leftarrow \nabla f(x_0)$
      %\STATE $\tilde{g} \leftarrow \nabla f(\tilde{x}_0)$ \\
      \FOR{$t = 0,1,\ldots,m-1$} 
	  %\STATE $\alpha_t \leftarrow \dfrac{(x_t^Tx_0)^{p-1}}{\|x_0\|^{2(p-1)}}$ \\
	  %\STATE $\alpha_t \leftarrow (x_t^Tx_0)^{p-1}$ \\
	  \STATE $\alpha_t \leftarrow {|x_t^Tx_0|^{p-1}} / {\|x_0 \|^{2(p-1)}}$ \\
	  \STATE sample $S_{t} \subset [n]$ of size $s$ uniformly at random \\
      %\STATE $g_t \leftarrow \dfrac{1}{|S_t|}\sum \limits_{l \in S_t} \nabla f_l(x_t) - \dfrac{\alpha_t}{|S_t|}\sum \limits_{l \in S_t} \nabla f_l(x_0) + \alpha_t \tilde{g}$ \\
      \STATE $\textstyle g_t \leftarrow \alpha_t \tilde{g}_T + s^{-1} \sum_{l \in S_t} \left[ \nabla f_l(x_t) - {\alpha_t} \nabla f_l(x_0) \right]$ \\
      \STATE $x_{t+1} \leftarrow (1-\eta) x_{t} + \eta {g_t}/{\| x_t \|^{p-2}}$ 
      \ENDFOR
      \STATE $\tilde{x}_{T+1} \leftarrow {x_m}$
      %\STATE $\tilde{x}_{s+1} \leftarrow \dfrac{x_m}{\|x_m\|_2}$
   \ENDFOR
\end{algorithmic}
\end{algorithm}

%In the additive scale invariant case, the algorithm remains the same except that we set $p=0$. The analyses remains the same since we only use the property $\nabla^2f(x)x = (p-1)\nabla f(x)$ for multiplicity while this expression for additive reads $\nabla^2f(x)x = -\nabla f(x)$. These expressions are provided in  \citep[Proposition 3]{kim2019scale}. Thus $p=0$ in the multiplicative case yields the additive case. 

%---------------------------------------------%
\section{Convergence Analysis}
\label{sec:convergence-analysis}
For the analysis of the algorithm, we assume that every $f_i$ is twice continuously differentiable on an open set containing $\partial \mathcal{B}_{d} \triangleq \{ y \in \mathbb{R}^d : \| y \| = 1 \}$. 

Let $x^*$ be a local optimal solution to \eqref{prob:main}. Then, by the first-order sufficient condition for optimality, there exists some $\lambda^*$ satisfying $\nabla f(x^*) = {\lambda^*} x^*$. Let $(\lambda_i,v_i)$ be an eigen-pair of $\nabla^2 f(x^*)$ and $\sigma = \| \nabla^2 f(x^*) \|$. Due to the eigenvector property \cite[Proposition 4]{kim2019scale} of the scale invariant problem, $x^*$ is an eigenvector of $\nabla^2 f(x^*)$. Without loss of generality, let $x^* = v_1$. 
Moreover, we assume that the local optimal $x^*$ statisfies $\lambda^* > \overline{\lambda} = {\textstyle \max_{2 \leq i \leq d}} |\lambda_i|$ as assumed in \cite[Theorem 7]
{kim2019scale}. 

Let $H_j$ be the Hessian of $\nabla_j f$ and 
$F_i(y^1,\cdots,y^d) = (\lambda^* - \lambda_1) \mathds{1}_{i = 1} I + \sum_{j=1}^d v_{ij} H_j(y^j)$.
Let $G_{S}({y}^1,\cdots,{y}^d)$ be the matrix such that $\nabla \nabla_j g_{S} ( {y}^j )^T$ is the $j^{th}$ row of $G_{S}({y}^1,\cdots,{y}^d)$ where $g_{S} = f_{S} - f$.

Next, we introduce some constants that are used to derive bounds in the analysis. Let $M = \max (M_1,M_2)$ where
\begin{equation} \label{def:M}
\begin{aligned}
    M_1 &= \hspace{-4mm} \max_{\substack{x \in \partial \mathcal{B}_d, \, y^1,\cdots,y^d \in \mathcal{B}_{d}}} \hspace{-2mm} \sqrt{\textstyle \sum_{i=1}^d ( x^T F_i(y^1, \cdots, y^d) x )^2}, \\
    M_2 &= \hspace{-4mm} \max_{\substack{x, z \in \partial \mathcal{B}_d, \, y^1,\cdots,y^d \in \mathcal{B}_{d}}} \hspace{-2mm} \big| \textstyle \sum_{i=1}^d z_i x^T H_i(y^1, \cdots, y^d) x \big|
\end{aligned}
\end{equation}
and $\mathcal{B}_{d} \triangleq \{ y \in \mathbb{R}^d : \| y \| \leq 1 \}$. These constants measure local smoothness of the objective function $f$ near $x^*$. Let $B_{s}$ be the set of all mini-batch samples $S \subset [n]$ of size $s=|S|$. We define quantities $K$ and $L$ as
\begin{align} \label{def:K-L}
\begin{aligned}
    K &= \hspace{-2mm} \max_{{y}^1,\cdots,{y}^d \in \mathcal{B}_{d}} \hspace{-2mm} E_{S \sim B_{s}} \, [\| G_{S}({y}^1,\cdots,{y}^d) \|^2], \\
    L &= \hspace{-4mm} \max_{\substack{S \in B_{s}, \, {y}^1,\cdots,{y}^d \in \mathcal{B}_{d} }} \, {\| G_{S}({y}^1,\cdots,{y}^d) \|^2}
\end{aligned}
\end{align}
and let $L_0$ be an upper bound of $L$ which can be obtained by setting $s=1$ (an easy calculation establish this). These constants measure deviation of $f_S$ from its mean $f$ with respect to stochastic sample $S$ of size $s$. $K$ measures the mean squared deviation (variance) of $f_S$ and $L$ is concerned with the maximum squared deviation of $f_S$ from $f$. As the batch size $s$ is increasing, both $K$ and $L$ are decreasing, and both of them become zero when $s=n$. While $K$ decreases as a factor of $1/s$, $L$ is a non-trivial function of $s$. Therefore, if some $f_i$ is extremely irregular (i.e. $|f_i-f|$ has an extremely large value around the solution), we would have to use a batch size close to $n$ to ensure that $L$ is smaller than some level as illustrated in Lemma~\ref{lemma:y0-St-necessary-condition}.

%{\color{orange} Let us emphasize that $K$ and $L$ are quantities depending on the batch size $s$, and measure irregularity of the worst-case stochastic gradient compared to the full-batch gradient. In general, $K$ and $L$ are non-trivial functions of $s$, and are exactly $0$ when $s = n$. When $f_i$'s are extremely irregular (i.e. $(f_i-f)$'s have extremely large variance around the solution), one has to use the full-batch.}

%{\color{orange} When $s = |S|$ is large, then $K$ and $L$ is small in general. It might be hard to quantify $K$ and $L$ in terms of $s$, since for PCA problem $G_S = s^{-1} \sum_{i \in S} x_i x_i^T - n^{-1}\sum_{i=1}^n x_ix_i^T$.}

Now, we present the convergence analysis for S-SCI-PI. We first analyze one-step inner iteration which computes $x_{t+1}$ from $x_t$. Let $\alpha(\eta) = 1-\eta+\eta \lambda^*$, $\beta (\eta) = 1-\eta+\eta \bar{\lambda}$, $y_k = x_k / \|x_k \|$ and $\Delta_t = 1 - y_t^Tx^*$. Since the optimality gap is expressed as $\sum_{i=2}^d (x_t^Tv_k)^2/(x_t^Tv_1)^2$, it is important to analyze how $x_t^Tv_k$ changes after each iteration. The following lemma provides an expression of $x_{t+1}^Tv_k$ as the sum of three components.
\begin{lemma}
\label{lemma:xt-vk-eq}
For $1 \leq k \leq d$ and any $t$, if $x_t^Tx_0 \geq 0$, then we have
\begin{align*}
    x_{t+1}^T v_k & = \left( 1-\eta+ \eta (\lambda_k + (\lambda^* - \lambda_1) \mathds{1}_{k = 1}) \right) x_t^Tv_k \\ 
    & +  \frac{1}{2} \eta \| x_t \| (y_t - x^*)^T F_k(\hat{y}_t^1,\cdots,\hat{y}_t^d) (y_t - x^*) \\
    & + \eta \left( G_{S_t}(\bar{y}_t^1,\cdots,\bar{y}_t^d) \big( x_t-(x_t^Ty_0)y_0  \big) \right)^T v_k
\end{align*}
for some $\hat{y}^1, \cdots, \hat{y}^d, \bar{y}^1, \cdots, \bar{y}^d \in \mathcal{B}_d$.
\end{lemma}
In Lemma~\ref{lemma:xt-vk-eq}, the first term represents the growth of $x_t^Tv_k$. The multiplicative factor is $1-\eta+\eta \lambda^*$ if $k=1$ and $1-\eta+\eta \lambda_k$ otherwise.
The second component is attributed to the difference of the Hessians at $x_t$ and $x^*$. As $x_t$ closes on $x^*$, this term goes to zero. The last term is stochastic error. The stochastic error is affected by the batch size $s$ and how closely $x_t$ is aligned with $x_0$ at which we compute the full gradient. 

The following lemma provides a condition on $\eta$, $L$, $M$ and $x_0$ to ensure that $y_t^Tx^*$ is not smaller than $y_0^Tx^*$ for every stochastic realization. 
\begin{lemma}
\label{lemma:y0-St-necessary-condition}
For any positive integer $m$, if the step size $\eta$, $s$ and $x_0$ are chosen to satisfy
\begin{align} 
    \label{Delta-0-condition}
   \Delta_0 {}& \leq \min \left \{ 1-\frac{1}{\sqrt{2}}, \, \frac{(\lambda^* - \bar{\lambda})^2}{4(M+2\sqrt{L_0})^2} \right \}
\end{align}
and either one of the following condition holds:
\\
\noindent\begin{minipage}{.5\linewidth}
\begin{equation} 
\label{L-condition}
    L \leq \frac{(\lambda^* - \bar{\lambda} - 2M \sqrt{\Delta_0})^2}{32},
\end{equation}
\end{minipage}
\begin{minipage}{.47\linewidth}
\begin{equation} 
\label{eta-condition} 
    \eta \leq 1 / \max (1, \nu_1, \nu_2, \nu_3),
\end{equation}
\end{minipage}
where
\begin{subequations}
\label{def:nu}
\begin{align}
    \nu_1 &= 1 - \lambda^* + 2 \theta_1 m \sqrt{2\Delta_0}, \label{def:nu-1}\\
    \nu_2 &= m\lambda^*+1-(m+1)(\bar{\lambda}+M\sqrt{\Delta_0}), \label{def:nu-2} \\
    \nu_3 &= \frac{128 L \theta_1 {\lambda^*}  m^2}{\theta_2^2 \bar{\lambda} \Delta_0 \sqrt{\Delta_0}} + 1 - \big( \bar{\lambda} + M \sqrt{\Delta_0} \big), \label{def:nu-3}
\end{align}
\end{subequations}
and
\begin{subequations}
\label{def:theta}
\begin{align}
    \theta_1 &= \lambda^* + \sigma + M \sqrt{\frac{\Delta_0}{2}} + 2\sqrt{L}, \label{def:theta-1} \\
    \theta_2 &= \lambda^* - \bar{\lambda} - 2 \sqrt{\Delta_0} \big( M + 2\sqrt{L} \big), \label{def:theta-2}
\end{align}
\end{subequations}
then we have $x_t^Tx_0 \geq 0$ and $\Delta_t \leq \Delta_0$ for all $0 \leq t \leq m$.
\end{lemma}
Note that $\lambda^*-\bar{\lambda}$ can be understood as a generalized eigen-gap at the solution, which specializes to $\lambda_1 - \lambda_2$ for the PCA problem. Since $L$ and $\Delta_0$ are decreasing functions of the batch size $s$ and the dot product $y_0^Tx^*$, given that $\Delta_0$ is moderately small, we can satisfy conditions \eqref{L-condition} or \eqref{eta-condition} by increasing the batch size $s$ or decreasing the step size $\eta$, respectively.
Conditioning on $x_t$, the next lemma derives expectation bounds for several quantities involving $(x_{t+1}^Tv_k)^2$ and norms.
\begin{lemma}
\label{lemma:recurrence-iteration}
For any positive integer $m$, if $\eta$, $s$ and $x_0$ satisfy \eqref{Delta-0-condition}, \eqref{L-condition} (or \eqref{eta-condition}) and
\begin{equation}
\label{eta-condition-1}
    \eta \leq 1/\max (1,1-\lambda^*+\sqrt{2} M \Delta_0 )
\end{equation}
then we have
\begin{align*}
    E[ \|x_{t+1}\|^2 | x_t] &\leq \big[ \big( \alpha(\eta) + \eta {M} \Delta_t \big)^2 + \eta^2 K \big] \|x_t\|^2, \\
    E \big[ \sum_{k=2}^d (x_{t+1}^T v_k)^2 | x_t \big] &\leq
    \big( \beta(\eta) + \eta M \sqrt{\Delta_t}  \big)^2  \sum_{k=2}^d (x_{t}^T v_k)^2 \\ 
    & + 8 \eta^2 K \|x_t\|^2  \sum_{k=2}^d (y_0^Tv_k)^2, \\
    E[ (x_{t+1}^Tv_1)^2 | x_t ] & \geq \bigg[ \alpha(\eta) - \frac{\eta M \Delta_t}{1-\Delta_t} \bigg]^2 (x_{t}^T v_1)^2
    %& E[ (x_{t+1}^Tv_1)^2 | x_t ] \leq \Big[ \alpha(\eta) + \frac{\eta M \Delta_t}{1-\Delta_t} \Big]^2 (x_{t}^T v_1)^2 \\
    %&\quad + 4\eta^2K \Big[ \sum_{k=2}^d (x_{t}^T v_k)^2 + \|x_t\|^2 \sum_{k=2}^d (y_0^Tv_k)^2 \Big].
\end{align*}
for any $0 \leq t \leq m$.
\end{lemma}
Using recursion on the one-step bound in Lemma~\ref{lemma:recurrence-iteration}, we derive bounds of $E[\sum_{k=2}^d (x_t^Tv_k)^2]$ and $E[(x_t^Tv_1)^2]$ as functions of $E[\sum_{k=2}^d (x_0^Tv_k)^2]$ and $E[(x_0^Tv_1)^2]$ below.
\begin{lemma}
\label{lemma:recurrence-epoch}
For any positive integer $m$, if $\eta$, $s$ and $x_0$ satisfy \eqref{Delta-0-condition}, \eqref{L-condition} (or \eqref{eta-condition}), \eqref{eta-condition-1} and
\begin{equation}
    \label{eta-condition-2}
    \eta \leq \max 1/(1, 1-\lambda^* - M \sqrt{\Delta_0} + \sqrt{Km}),
\end{equation}
then we have
\begin{align*}
    E \big[ \sum_{k=2}^d (x_{t}^T v_k)^2 \big] & \leq
    E \big[ \sum_{k=2}^d (x_0^Tv_k)^2 \big]
    \big[ \big( \beta(\eta) + \eta M \sqrt{\Delta_0}  \big)^{2t} \\
    & + 16 \eta^2 K t \big( \alpha(\eta) + \eta {M} \sqrt{\Delta_0} \big)^{2(t-1)} \big], \\
    E[ (x_{t}^Tv_1)^2] &\geq \bigg[ \alpha(\eta) - \frac{\eta M \Delta_0}{1-\Delta_0} \bigg]^{2t} E[(x_{0}^T v_1)^2].
\end{align*}
\end{lemma}
The inequalities in Lemma~\ref{lemma:recurrence-epoch} are important since they yield a bound on the optimality gap which is expressed as $\delta_t \triangleq E [ \sum_{k=2}^d (x_{t}^T v_k)^2]/E[ (x_{t}^Tv_1)^2]$. In the next lemma, we show that under some conditions on $\eta, m, s$ and $x_0$, the optimality gap decreases at least by $1-\rho$ after each outer iteration. 

\begin{lemma}
\label{lemma:expectation-bound-epoch}
For any positive integer $m$, if $\eta$, $s$ and $x_0$ satisfy
\eqref{Delta-0-condition}, \eqref{L-condition} (or \eqref{eta-condition}) and
\begin{align}
    \label{eta-condition-3}
    \eta \leq 1/\max(1, \nu_4, \nu_5)
\end{align}
where
\begin{subequations}
\begin{align}
    \label{def:nu-4}
    \nu_4 &= 1-\lambda^* - M \sqrt{\Delta_0} \\
    & + \max 
    \left( \sqrt{Km}, \, \frac{64 K}{\lambda^* - \bar{\lambda} - 2M\sqrt{\Delta_0}} \right) \nonumber \\ 
    \label{def:nu-5}
    \nu_5 &= 1-\lambda^* + M \sqrt{\Delta_0} \\
    & + \max 
    \left( 2m \big(\lambda^* - \bar{\lambda} - 2M\sqrt{\Delta_0} \big), \, \frac{4m M \sqrt{\Delta_0}}{\log 2}  \right), \nonumber
\end{align}
\end{subequations}
then we have $\delta_m \leq (1- \rho) \cdot \delta_0$
where
\begin{equation}
    \label{def:rho}
    0 < \rho = \frac{\eta m \big( \lambda^* - \bar{\lambda} - 2M\sqrt{\Delta_0} \big)}{2 \big( 1- \eta + \eta (\lambda^* - M \sqrt{\Delta_0}) \big)} < 1.
\end{equation}
\end{lemma}
Finally, we analyze the entire algorithm. Let $\tilde{\Delta}_0 = 1- \tilde{x}_0^Tx^*$ and $\tilde{\delta}_s = E[ \sum_{k=2}^d (\tilde{x}_s^Tv_k)^2]/E[(\tilde{x}_s^Tv_1)^2]$.
By repeatedly applying Lemma~\ref{lemma:expectation-bound-epoch}, the following theorem states that $\tilde{\delta}_s$ decreases at a liner rate under some conditions on $\eta$, $m$, $s$ and $\tilde{x}_0$.
\begin{theorem}
\label{thm:convergence}
For any positive integer $m$, if $\eta$, $s$ and $\tilde{x}_0$ satisfy
\eqref{Delta-0-condition}, \eqref{L-condition} (or \eqref{eta-condition}) and \eqref{eta-condition-3} with $\Delta_0 = \tilde{\Delta}_0$, 
then for any $\epsilon > 0$, after $\tau = \lceil (1/\rho) \log( \tilde{\delta}_0/\epsilon) \rceil$ epochs of S-SCI-PI (Algorithm~\ref{alg:Stochastic-SCI-PI}), we have $\tilde{\delta}_\tau \leq \epsilon$.
\end{theorem}
Theorem~\ref{thm:convergence} states that for any epoch length $m$, if $\widetilde{x}_0$ is moderately close to $x^*$ and the step size $\eta$ and the batch size $s$ satisfies certain conditions, the optimality gap vanishes at an exponential rate. If there are few irregular $f_i$ and the cost of sampling is low, we can satisfy \eqref{L-condition} by making $L$ small. In this case, $\eta$ can take a large value and we are able to obtain rapid convergence. On the other hand, if there are many irregular data samples and sampling is expensive, we may not be able to satisfy \eqref{L-condition}. Nevertheless, we can always ensure linear convergence of Algorithm~\ref{alg:Stochastic-SCI-PI} by choosing a small enough step size $\eta$ (conditions \eqref{eta-condition}, \eqref{eta-condition-1}, \eqref{eta-condition-2}, \eqref{eta-condition-3}) as in \cite{shamir2015stochastic}.
%---------------------------------------------%
\section{Application: KL-divergence NMF}
\label{sec:KL-NMF}
%---------------------------------------------%
Let $V \in \mathbb{R}_+^{N \times M}$ be a given non-negative matrix, which we want to compress into the product of $W \in \mathbb{R}_+^{N \times K}$ and $H \in \mathbb{R}_+^{K \times M}$. The KL-NMF problem is defined as
\begin{equation*} %\label{eqn:kldiv_nmf}
\begin{aligned}
& \min & & \sum_{i,j} \left[V_{ij} \log \frac{V_{ij}}{(WH)_{ij}} - V_{ij} + (WH)_{ij} \right] \\
& \text{subject to} & & W_{ik}\geq 0, \,\, H_{kj} \geq 0, \,\, \forall i,j,k.
\end{aligned}
\end{equation*}
The above objective function is called the (generalized) Kullback-Leibler divergence $D_{KL}(V\| WH)$. Let $H_j$ be the $j$-th column of $H$. Note that the objective function $D_{KL}(V\| WH)$ is separable in $H_1, \cdots, H_M$ and thus
\baa \label{prob:kldivsubproblem}
\begin{aligned}
H_j^{\rm new} = \,\, & \argmax & & {\textstyle \sum_i} \left[V_{ij} \log (WH_j)_i - (WH_j)_i\right] \\
\,\, & \text{subject to } & & H_j \geq 0
\end{aligned}
\eaa
serves as the $j$-th subproblem.
%---------------------------------------------%
\begin{lemma}{\cite[Lemma 9]{kim2019scale}} \label{lem:klequivalence}
The $j$-th KL-NMF subproblem \eqref{prob:kldivsubproblem} is equivalent to:
\baa \label{prob:mixsqpsubproblem}
\begin{aligned}
X_j^{\rm new} = \,\, & \argmax & & {\textstyle\sum_{i=1}^N} V_{ij} \log (LX_j)_i \\
& \textrm{\rm subject to} & & X_j \in \mathcal{S}^d
\end{aligned}
\eaa
where $L_{ik} = W_{ik} / (\sum_{i'} W_{{i'}k})$ and $\mathcal{S}^d \triangleq \{ x : \textstyle \sum_{k=1}^d x_i = 1, \, x \succeq 0 \}$. The original solution $H_{j}^{\rm new}$ can be recovered via
\baa
    H_{kj}^{\rm new} = \frac{\sum_i V_{ij}}{\sum_i W_{ik}} \, X_{jk}^{\rm new}, \quad k \in [d].
\eaa
\end{lemma}
%---------------------------------------------%
By re-parameterizing $X_j$ by $Y_j^2$, we can convert \eqref{prob:mixsqpsubproblem} into a finite-sum scale invariant problem \eqref{prob:main}.
%Lemma~\ref{lem:klequivalence} implies that the KL divergence NMF subproblem for $H$, namely
Therefore, the KL divergence NMF subproblem for $H$, namely
\ba
\begin{aligned}
H^{\rm new} \,\, = & \argmax & & [V_{ij} \log (WH)_{ij} - (WH)_{ij}] \\
& \text{subject to} & & {H \geq 0}
\end{aligned}
\ea
can be solved by applying S-SCI-PI for each column of $H$.
%can be solved by S-SCI-PI after we reformulate problem \eqref{prob:kldivsubproblem} into \eqref{prob:mixsqpsubproblem} for $j \in [M]$.

The vanilla stochastic updates $H$ by sampling a mini-batch $S \leftarrow \textrm{sample}(N,s)$ and running
\baa \label{eqn:vanilla-sscipi}
\begin{aligned}
H_{kj}^{\rm new} & \leftarrow H_{kj} \left[(1-\eta) + \eta \sum_{i \in S} \frac{L_{ik} V_{ij}}{(LH)_{ij}}\right]^2 \,\, \forall k, \, j, \\
H^{\rm new} & \leftarrow \textrm{column-rescale}(H^{\rm new}),
\end{aligned}
\eaa
where sample$(N,s)$ is sampling $s$ elements from $[N]$ with or without replacement, and column-rescale$(X)$ is rescaling the columns of $X$ to have sum $1$. The update for $W$ is similar due to
$D_{KL}(V \| WH) = D_{KL} (V^T \| H^T W^T)$. Algorithm~\ref{alg:Stochastic-SCI-PI} (S-SCI-PI) can be understood as an SVRG version of \eqref{eqn:vanilla-sscipi}.

Our final remark is that we solve the $j$-th KL-NMF problem for $j \in [M]$ simultaneously as a single optimization problem. Let $X = [X_1,\cdots,X_M]$ be the concatenation of the $M$ column vectors $X_1,\cdots,X_M \in \mathbb{R}^K$. Lemma~\ref{lem:klequivalence} states that in the exact alternating minimization algorithm, the update of $H$ amounts to solving
\baa \label{prob:combined-klnmf-subprob}
\begin{aligned}
& \min & & {\textstyle \sum_{i = 1}^{NM} \textrm{vec}(V)_i \log [(I_M \otimes L) \textrm{vec}(X)]_i} \\
& \textrm{subject to} & & \textrm{vec}(X_j) \in \mathcal{S}^K,\ \ j \in [M]
\end{aligned}
\eaa
where $\textrm{vec}(X) = (X_{11},\cdots,X_{K1},\cdots,X_{1M},\cdots,X_{KM})$ is a vectorization of $X \in \mathbb{R}^{K \times M}$, $\textrm{vec}(V) \in \mathbb{R}^{KM}$ is defined similarly and $I_M \otimes L = \textrm{kron}(I_M,L) \in \mathbb{R}^{NM \times KM}$ is the Kronecker product of $I_M$ and $L$.

% i.e.
% \ba
% I_M \otimes L = \textrm{kron}(I_M,L) = \left[\begin{array}{cccc}
% L & \cdots & 0 \\
% \vdots & \ddots & \vdots \\
% 0 & \cdots & L \end{array} \right]
% \ea
% \begin{lemma} \label{lem:vectorized-prob}
% The solution to \eqref{prob:combined-klnmf-subprob} can be easily recovered once we solve
% \baa \label{eqn:vectorized-prob}
% \textrm{minimize} {}& \quad \sum_{i = 1}^{NM} \textrm{\rm vec}(V)_i \log [(I_M \otimes L) \textrm{\rm vec}(X)]_i, \\
% \textrm{subject to} {}& \quad \textrm{\rm vec}(X) \in \mathcal{S}^{KM}.
% \eaa
% and rescale columns of $X$ to have sum $1$ so that $X_j \in \mathcal{S}^K$.
% \end{lemma}
This allows us to exploit fast matrix multiplication routines (i.e. efficient matrix computation library such as OpenBLAS or intel MKL) in solving the aggregated problem \eqref{prob:combined-klnmf-subprob}, instead of solving the $j$-th subproblem sequentially for $j \in [M]$.
\iffalse
%---------------------------------------------%
\subsection{Algorithm}
%---------------------------------------------%
The gradient of \eqref{prob:combined-klnmf-subprob} is
\ba
L^T [V \oslash (LX)]
\ea
%---------------------------------------------%
\subsection{Adaptive Sampling}
%---------------------------------------------%
\ba
f(x)  = \sum_{i=1}^{n} v_i f_i(x),\quad \nabla f(x) = \sum_{i=1}^n v_i \nabla f_i(x)
\ea
We will write $S \sim v^s$ when $|S| = s$ and each element in $S$ is drawn from $\{1,\cdots,n\}$ with probability $(v_1,v_2,\cdots,v_n) / \sum_{i=1}^n v_i$. In other words, $S$ is a random sample of size $s$ drawn from the Categorical$(v_1,\cdots,v_n)$ distribution on $\{1,\cdots,n\}$. Then
\ba
\mathbb{E}_{S \sim v^s} \left[{\sum_{l \in S}} \nabla f_l(x) \right] = \nabla f(x) 
\ea
\fi
%---------------------------------------------%
\subsection{Related Algorithms}
\label{subsec:relatedmethod}
%---------------------------------------------%
Let $Z = WH$ henceforth. We omit the update of $W$ since it can be derived similarly.

\textit{Multiplicative Update (MU) / Expectation Maximization (EM) \citep{lee2001algorithms}}:
MU updates all $H_{kj}$'s simultaneously by
\ba
H_{kj}^{\rm new} = H_{kj} \frac{\sum_{i} W_{ik} V_{ij} / Z_{ij}}{\sum_i W_{ik}}
\ea
for all $k$ and $j$. \iffalse, and then updates all $W_{ik}$'s by
\ba
W_{ik}^{\rm new} = W_{ik} \frac{\sum_{j} H_{kj} V_{ij} / Z_{ij}}{\sum_j H_{kj}}
\ea
for all $i$ and $k$.
\fi
Let us emphasize that the MU update is identical to the standard EM algorithm for the estimation of mixture proportions.

\textit{Cyclic/Stochastic Coordinate Descent (CCD/SCD)} \citep{hsieh2011fast,muzzarelli2019rank}:
For all $j$ and $k$, CCD/CSD runs coordinate-wise updates of $H$
\ba
H_{kj}^{\rm new} {}& = \max\left\{0, H_{kj} - \frac{\sum_i W_{ik} (1 - V_{ij} / Z_{ij})}{\sum_i V_{ij} W_{ik}^2 / Z_{ij}^2 } \right\}
\ea
sequentially in a pre-fixed cyclic order.\iffalse Then for $i=1,\cdots,N$ and $k = 1,\cdots,K$, run coordinate-wise updates of $H$
\ba
W_{ik}^{\rm new} {}& = \max\left\{0, W_{ik} - \frac{\sum_j H_{kj} (1 - V_{ij} / Z_{ij})}{\sum_j V_{ij} H_{kj}^2 / Z_{ij}^2 }  \right\}
\ea
sequentially in a pre-fixed cyclic order.\fi

\textit{Projected Gradient Descent (PGD)} \citep{lin2007projected}:
Given element-wise step sizes $\alpha_{kj}$'s, PGD updates all $H_{kj}$'s simulataneously via
\ba
%\textstyle 
H_{kj}^{\rm new} = \max\left\{0,  H_{kj} - \alpha_{kj} \sum_{i} W_{ik} \left( 1 - \frac{V_{ij}}{Z_{ij}} \right) \right\}.
\ea
\iffalse
and then updates all $W_{ik}$'s via
\ba
W_{ik}^{\rm new} = \max \left\{0,  H_{ik} - \alpha_{ik} \left(\sum_{i} H_{kj} (1 - V_{ij} / Z_{ij}) \right) \right\}
\ea
\fi
Note that Multiplicative Update (MU) is a special case of PGD when $\alpha_{kj} = H_{kj} / (\sum_i W_{ik})$, which does not require projection onto the non-negative orthant. Also, CCD updates $H_{kj}$ one at a time with a coordinate-wise optimal step size $\alpha_{kj} = 1 / \sum_i (V_iW_{ik}^2 / Z_{ij}^2)$. By contrast, PGD uses a single step size $\alpha_j = \alpha_{1j} = \cdots = \alpha_{Kj}$ for each column $j$ for fast line searches.

%The proposed S-SCI-PI algorithm (\ref{eqn:vanilla-sscipi}) is most similar to MU since $W$ and $H$ are updated multiplicatively as well. As a special case of S-SCI-PI, the full-batch S-SCI-PI is denoted by F-SCI-PI.

Let us highlight that S-SCI-PI and all the comparison methods belong to the family of alternating minimization algorithms, which update $H$ given $W$ and then update $W$ given $H$ iteratively.
\begin{figure*}[ht]
    \centering
    \includegraphics[width = 5.75in]{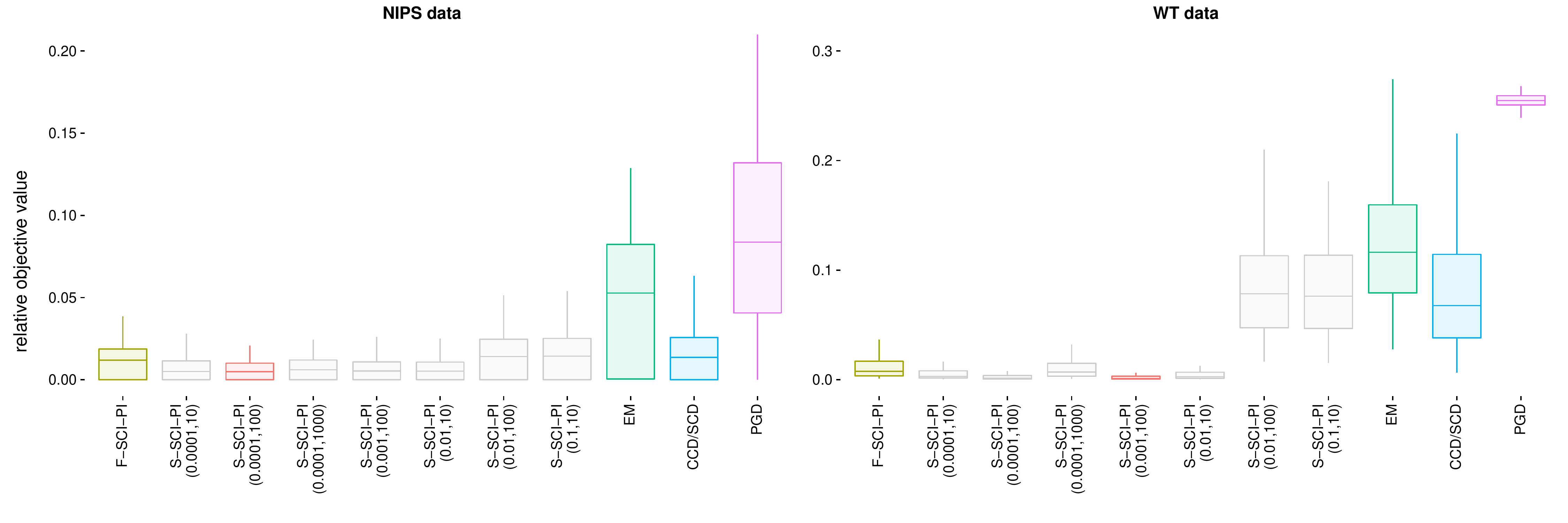}
    \caption{Boxplots of the relative errors after 30 seconds across 10 independent replicates}
    %\caption{Boxplots showing variability of each method across 10 independent replicates.}
    \label{fig:variance}
\end{figure*}
%---------------------------------------------%
\subsection{Practical Considerations}
%\paragraph{Various Sampling Scheme}
\label{sec:practical-considerations}
%---------------------------------------------%
In this part, we compare several sampling schemes for the update of $H$. Since the sampling scheme for the update of $W$ can be similarly discussed, we omit it.

\textit{Vector-wise Sampling: }We construct $V_S \in \mathbb{R}_+^{s \times M}$ and $W_S \in \mathbb{R}_+^{s \times D}$ by sampling rows of $V \in \mathbb{R}_+^{N \times M}$ and $W \in \mathbb{R}_+^{N \times D}$ uniformly at random, respectively. The stochastic gradient reads
\ba
\nabla_S^{\rm row} f(H) = \frac{n}{s} W_S^T [V_S \oslash (W_S H)].
\ea
For a dense data matrix $V$, we prefer to use this vector-wise (or row-wise) sampling scheme for the update of $H$, since it allows us to exploit fast matrix multiplication libraries.

\textit{Element-wise Sampling: } We vectorize the problem by introducing the element-wise iterator $i = (i_1,i_2) \in [N] \times [M] = [NM]$. This yields
\ba
f(H) = \sum_{i \in \mathcal{I}} V_{i_1, i_2} \log \sum_{k=1}^D W_{i_1,k} H_{k,i_2}
\ea
where $\mathcal{I}$ is the subset of $[NM]$ such that $V_{i_1,i_2} \neq 0$ if and only if $i = (i_1,i_2) \in \mathcal{I}$. In other words, $\mathcal{I}$ is the index set of the nonzero elements in $V$.

We construct $S$ by sampling $s$ elements of $\mathcal{I}$ uniformly at random, and consider the stochastic gradient as
\ba
\nabla_S^{\rm elem} f(H) = \frac{|\mathcal{I}|}{s} \sum_{i \in S} \sum_k \frac{W_{i_1,k} V_{i_1, i_2}}{\sum_{k'=1}^D W_{i_1,k'} H_{k',i_2}} E_{k,i_2}
\ea
where $E_{k,i_2}$ is the standard basis matrix having $1$ at $(k,i_2)$-th entry and $0$ otherwise.

For a sparse data matrix $V$, we prefer this element-wise sampling scheme for $H$ over the row-wise sampling scheme, since each column has a different sparsity pattern.
% \paragraph{Adaptive Sampling: } For KL-NMF problem, each entries of $V$ has different importance due to the sparsity pattern of $V$ as well as the magnitude of rows, columns and elements of $V$. In words, one may prefer to sample ``keywords'' more often than ``common words'' while keeping the expectation of the stochastic gradient is equal to the full-batch gradient.

% For instance, one can sample larger entries more frequently. Also, one can sample $V_{ij}$'s more often when the $i$-th row vector $V_i$ and the $j$-th column vector $V_j$ of $V$ are more sparse. Finding a good adpative sampling scheme would be an interesting future direction to pursue.

\textit{Numerical Stability:}
The KL-NMF objective function and its gradient are unstable when entries of $V$ and $WH$ are close to $0$. 
%As reported in \citet{kim2019scale} and based on the experiments in Section~\ref{sec:experiment}, MU and the full-batch version of F-SCI-PI are numerically stable. 
For instance, if stochastic samples have many zero entries, stochastic variance-reduced gradients can have some negative elements, which may lead to numerical errors when computing objective values. 
In order to prevent this numerical issue, we use element-wise sampling and set a threshold to ensure that stochastic variance-reduced gradients are always non-negative.
%---------------------------------------------%
\section{Experiment}
\label{sec:experiment}
%---------------------------------------------%
\iffalse
\begin{figure}[ht]
    \centering
    \includegraphics[width = 3.5in]{figure3_for_paper.pdf}
    \caption{Caption}
    \label{fig:sim}
\end{figure}
\fi
\begin{figure*}[!t]
    \centering
    \iffalse
    \includegraphics[width = 5.25in]{figure1_for_paper.pdf}
    \vspace{-0.2in}
    \caption{Convergence plots (relative error vs. computation time) of the one-step alternating minimization algorithms on synthetic data sets. A few initial data points are removed for increased visibility.}
    \label{fig:sim}
    \fi
    \includegraphics[width = 5.5in]{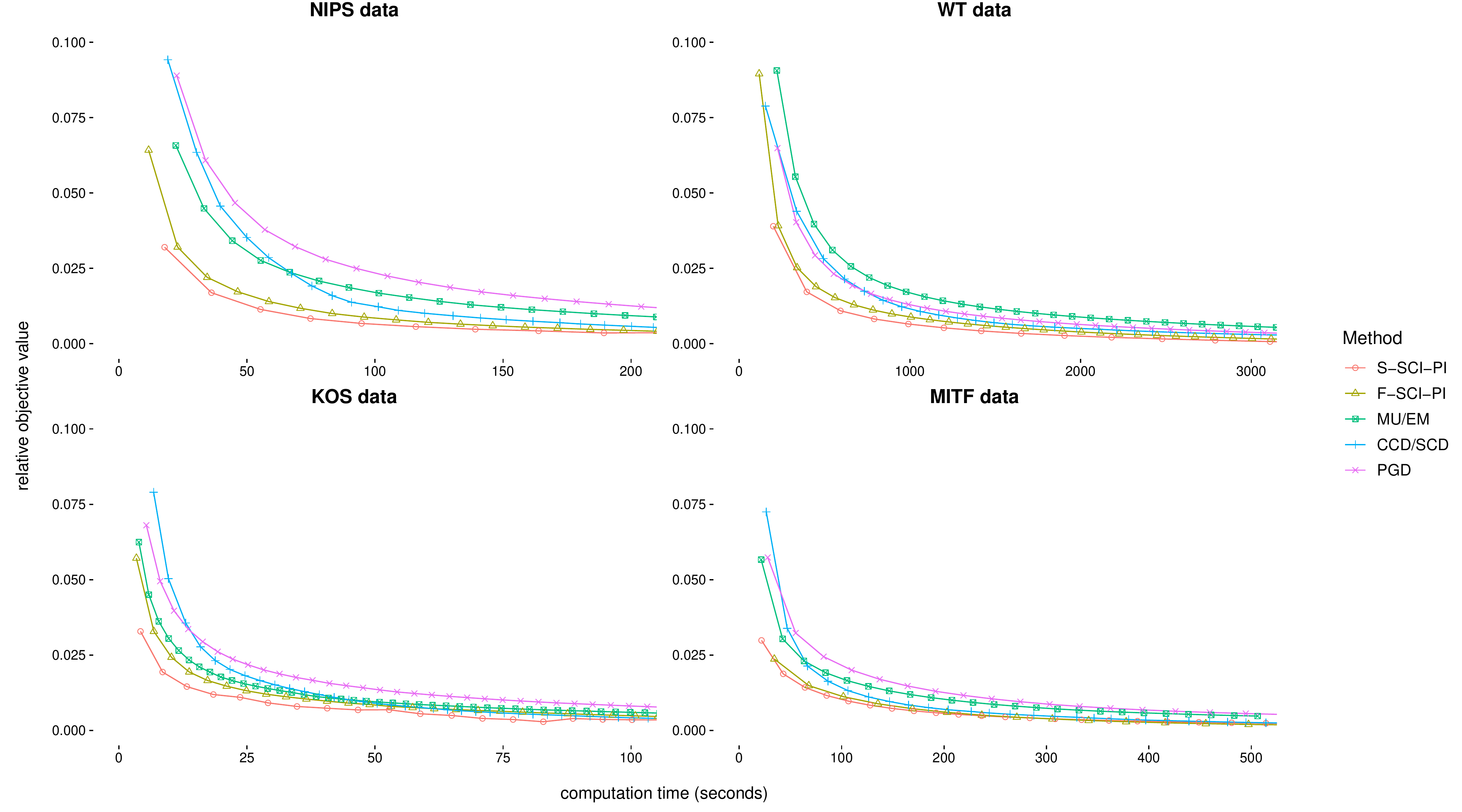}
    \caption{Convergence plots (relative error vs. computation time) of KL-NMF one-step alternating minimization algorithms on real data sets.}
    \label{fig:realnmf}
    \iffalse
    \vspace{0.2in}
    \includegraphics[width = 5.5in]{figure21_for_paper.pdf}
    \vspace{-0.2in}
    \caption{Convergence plots (relative error vs. iteration) of one-step alternating minimization on real data sets.}
    \label{fig:realnmfiter}
    \vspace{-0.1in}
\fi
\end{figure*}
We test the proposed algorithm S-SCI-PI on synthetic and real-world data sets. All experiments are implemented on a standard laptop (2.6 GHz Intel Core i7 processor and 16GB of RAM) using the C\texttt{++} programming language. We use 4 real data sets publicly available online and 3 synthetic data sets generated from Poisson distributions. The description is provided in Appendix~\ref{sec:data}. We set $K = 20$ features. All the reported values are averaged over 10 independent replicates started at different initial points, each of which is obtained by running 5 MU/EM steps on a Uniform(0,1) random matrix.
% In certain runs due to numerical errors the outcomes were peculiar and thus they were disregarded (but each observation has 5 normal runs). The benchmark algorithms are MU/EM, CCD/SCD, PGD, and F-SCI-PI. 
For S-SCI-PI, we perform grid search on the parameters by selecting the best parameters among different batch proportions $s/n \in \{0.0001,0.001,0.01,0.1, 1\}$, epoch lengths $m \in \{10, 100, 1000\}$ and step sizes $\eta \in \{0.01,0.1,1\}$.

\paragraph{Error metric}
We report the relative objective value (or the relative error) defined as
\ba
\Delta_t^{\rm rel} = |f(x^t) - f(x^*)| / |f(x^0) - f(x^*)|,
\ea
which always makes the reported values equal to one for all data sets at initial points. Note that this is simply an affine transformation of the standard error metric. This relative error is introduced to compare performance of the methods on different data sets.
\begin{figure*}[!t]
    \centering
    \includegraphics[width = 6in]{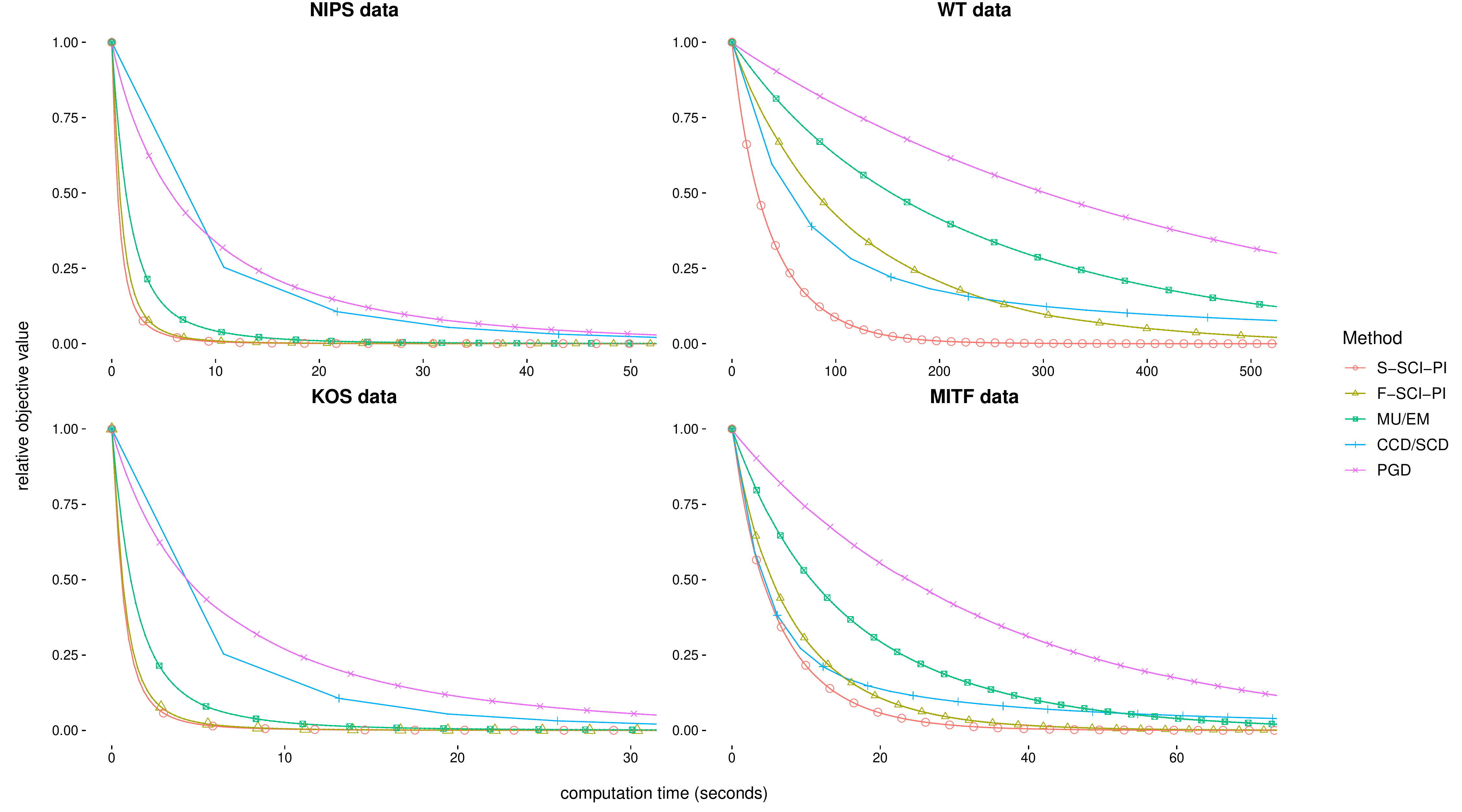}
    \caption{Convergence plots (relative error vs. computation time) for the KL-NMF subproblem.}
    \label{fig:exact}
    \iffalse
    \vspace{0.1in}
    \centering
    \includegraphics[width = 5.25in]{figure4_for_paper.pdf}
    \caption{Convergence plots (relative error vs. computation time) for the KL-NMF subproblem.}
    \label{fig:variance}
        \fi
    \vspace{-0.1in}
\end{figure*}

\iffalse
\begin{figure*}[ht]
    \centering
    \includegraphics[width = 6.25in]{figure1_for_paper.pdf}
    \vspace{-0.1in}
    \caption{Convergence plots (relative error vs. computation time) of the one-step alternating minimizations on the synthetic data sets. All the reported values are from every 50, 50, 100 and 100 iterations, respectively from the left to the right.}
    \label{fig:sim}
        \vspace{-0.1in}
\end{figure*}
\fi
%---------------------------------------------%
\paragraph{KL-NMF (one-step alternating minimization):}
%---------------------------------------------%
The one-step alternating minimization scheme is to update $H$ via a single iteration of each algorithm and then update $W$ similarly. A single iteration amounts to one outer loop iteration for S-SCI-PI, which involves $m$ stochastic updates.
We compare S-SCI-PI, F-SCI-PI, MU/EM, CCD/SCD and PGD.

For dense data sets (WT, MITF), we apply vector-wise sampling only on the columns (of dimension 19,200 and 2,429, respectively) since the other dimension is small (287 and 361, respectively). For sparse data sets (NIPS, KOS), the element-wise sampling scheme is applied to both dimensions, which turns out to be more effective.

%Figure~\ref{fig:sim} displays the results for the 4 synthetic data sets. By comparing the left two figures in Figure~\ref{fig:sim}, and the right figures, we conclude that S-SCI-PI performs much better than F-SCI-PI for dense matrices, but is the winner also for sparse matrices.

Figure~\ref{fig:realnmf} displays the relative errors with respect to the computation time for the 4 real data sets. Overall, S-SCI-PI with the chosen batch and epoch size improves the convergence over F-SCI-PI. However, S-SCI-PI does not outperform F-SCI-PI for the MITF data set, which has a relatively small number of columns (2,429). Also, both S-SCI-PI and F-SCI-PI exhibit much faster convergence than MU/EM. This clearly attests that S-SCI-PI is a reliable and practical option for the KL-NMF problem.

%Figure~\ref{fig:realnmfiter} displays the relative errors with respect to the outer loop iterations. S-SCI-PI takes longer steps per outer loop iteration than F-SCI-PI and MU/EM, at the expense of larger computational complexity. The batch size and the epoch length balance the trade-off between them. We also notice that F-SCI-PI and MU/EM have almost the same computation time but F-SCI-PI takes longer step than MU/EM.
%---------------------------------------------%
\paragraph{KL-NMF subproblem (exact alternating minimization):}
%---------------------------------------------%
The exact alternating minimization scheme is to update $H$ until it reaches the exact coordinate minimizer and then update $W$ similarly.
Instead of solving the entire KL-NMF problem, we solve a single KL-NMF subproblem to optimality and plot relative objective values over time to compare the speed of convergence.
%Thus we empirically study optimization performance of the algorithms on the KL-NMF subproblem, which is equivalent to the estimation of mixture proportions \citep{kim2019scale}. To this end, we first compare the 5 competitors.

Figure~\ref{fig:exact} displays the results for the real world data sets. It shows that S-SCI-PI is an overall winner solving the KL divergence subproblems and hence an efficient method for exact alternating minimization. However, it does not outperform F-SCI-PI significantly on the sparse NIPS and KOS data set. As reported in \cite{hsieh2011fast}, CCD/SCD is faster than MU/EM for the dense WT data set. However, our result on NIPS and KOS shows that CCD/SCD is much slower than S-SCI-PI mainly due to the expensive coordinate updates.

\paragraph{Robustness of S-SCI-PI:} Lastly, we compare the performance of S-SCI-PI for select choices of batch proportion $s/n$ and epoch length $m$. For each choice of $s/n$ and $m$, we find the best step size $\eta$ using the same grid search as above. For the NIPS and WT data sets, we run the algorithms for 30 seconds with $10$ independent replications and report boxplots of $10$ relative objective values in Figure~\ref{fig:variance}. In this figure, S-SCI-PI ($n/s$,$m$) stands for S-SCI-PI with batch proportion $n/s$ and epoch length $m$. The boxplots show that the performance of S-SCI-PI is robust to batch proportion $s/n$ and epoch length $m$ given that step size $\eta$ is appropriately selected.
% We compare the convergence of the methods. 
%We select two data sets (NIPS and WT) and report the relative objective values of each method, including S-SCI-PI with few selected parameters (different batch sizes, epoch lengths and step sizes). 
%The result is summarized in Figure~\ref{fig:variance}. 
%They show that a careful choice of the batch size $s$ and epoch length $m$ yields non-negligible improvements on the convergence of S-SCI-PI.
%Again, we confirm that the stochastic approach (S-SCI-PI) has a remarkable improvement over the full gradient approach (F-SCI-PI) for the dense WT data set.
Again, we emphasize that S-SCI-PI has a remarkable improvement over the full gradient approach (F-SCI-PI) for the dense WT data set.
%---------------------------------------------%
\section{Final Remarks}
%---------------------------------------------%
We introduce a stochastic variance-reduced algorithm (S-SCI-PI) to solve finite-sum scale invariant problems for the first time in the literature and provide its convergence analysis. Our analysis shows that under some conditions on the step size, epoch length, batch size and initial iterate, the algorithm achieves linear convergence in expectation. Using S-SCI-PI, we introduce a stochastic approach to solve the KL-NMF problem. The experimental results reveal that S-SCI-PI exhibits robust and superior performance over state-of-the-art methods.
\newpage
\bibliographystyle{plainnat}
\bibliography{main}

\begin{thebibliography}{31}
\providecommand{\natexlab}[1]{#1}
\providecommand{\url}[1]{\texttt{#1}}
\expandafter\ifx\csname urlstyle\endcsname\relax
  \providecommand{\doi}[1]{doi: #1}\else
  \providecommand{\doi}{doi: \begingroup \urlstyle{rm}\Url}\fi

\bibitem[Arora et~al.(2012)Arora, Cotter, Livescu, and
  Srebro]{arora2012stochastic}
Raman Arora, Andrew Cotter, Karen Livescu, and Nathan Srebro.
\newblock {Stochastic Optimization for PCA and PLS}.
\newblock In \emph{Annual Allerton Conference on Communication, Control, and
  Computing}, pages 861--868. IEEE, 2012.

\bibitem[Arora et~al.(2013)Arora, Cotter, and Srebro]{arora2013stochastic}
Raman Arora, Andy Cotter, and Nati Srebro.
\newblock {Stochastic Optimization of PCA with Capped MSG}.
\newblock In \emph{Advances in Neural Information Processing Systems}, pages
  1815--1823, 2013.

\bibitem[Boutsidis et~al.(2015)Boutsidis, Garber, Karnin, and
  Liberty]{boutsidis2015online}
Christos Boutsidis, Dan Garber, Zohar Karnin, and Edo Liberty.
\newblock {Online Principal Components Analysis}.
\newblock In \emph{Proceedings of the Twenty-Sixth Annual ACM-SIAM Symposium on
  Discrete Algorithms}, pages 887--901. Society for Industrial and Applied
  Mathematics, 2015.

\bibitem[Erdogdu et~al.(2018)Erdogdu, Ozdaglar, Parrilo, and
  Vanli]{erdogdu2018convergence}
Murat~A Erdogdu, Asuman Ozdaglar, Pablo~A Parrilo, and Nuri~Denizcan Vanli.
\newblock {Convergence Rate of Block-Coordinate Maximization Burer-Monteiro
  Method for Solving Large SDPs}.
\newblock \emph{arXiv preprint arXiv:1807.04428}, 2018.

\bibitem[F{\'e}votte and Idier(2011)]{fevotte2011algorithms}
C{\'e}dric F{\'e}votte and J{\'e}r{\^o}me Idier.
\newblock {Algorithms for Nonnegative Matrix Factorization with the
  $\beta$-divergence}.
\newblock \emph{Neural Computation}, 23\penalty0 (9):\penalty0 2421--2456,
  2011.

\bibitem[Garber et~al.(2015)Garber, Hazan, and Ma]{garber2015online}
Dan Garber, Elad Hazan, and Tengyu Ma.
\newblock {Online Learning of Eigenvectors}.
\newblock In \emph{International Conference on Machine Learning}, pages
  560--568, 2015.

\bibitem[Hardt and Price(2014)]{hardt2014noisy}
Moritz Hardt and Eric Price.
\newblock {The Noisy Power Method: A Meta Algorithm with Applications}.
\newblock In \emph{Advances in Neural Information Processing Systems}, pages
  2861--2869, 2014.

\bibitem[Hsieh and Dhillon(2011)]{hsieh2011fast}
Cho-Jui Hsieh and Inderjit~S Dhillon.
\newblock {Fast Coordinate Descent Methods with Variable Selection for
  Non-negative Matrix Factorization}.
\newblock In \emph{Proceedings of the 17th ACM SIGKDD international conference
  on Knowledge discovery and data mining}, pages 1064--1072, 2011.

\bibitem[Hyvarinen(1999)]{hyvarinen1999fast}
Aapo Hyvarinen.
\newblock {Fast ICA for Noisy Data using Gaussian Moments}.
\newblock In \emph{Proceedings of the 1999 IEEE International Symposium on
  Circuits and Systems VLSI}, volume~5, pages 57--61. IEEE, 1999.

\bibitem[Hyv{\"a}rinen et~al.(2004)Hyv{\"a}rinen, Karhunen, and
  Oja]{hyvarinen2004independent}
Aapo Hyv{\"a}rinen, Juha Karhunen, and Erkki Oja.
\newblock \emph{Independent Component Analysis}, volume~46.
\newblock John Wiley \& Sons, 2004.

\bibitem[Johnson and Zhang(2013)]{johnson2013accelerating}
Rie Johnson and Tong Zhang.
\newblock {Accelerating Stochastic Gradient Descent using Predictive Variance
  Reduction}.
\newblock In \emph{Advances in Neural Information Processing Systems}, pages
  315--323, 2013.

\bibitem[Jolliffe(2002)]{jolliffe2002principal}
Ian Jolliffe.
\newblock \emph{Principal \uppercase{C}omponent \uppercase{A}nalysis}.
\newblock Wiley Online Library, 2002.

\bibitem[Kasai(2018)]{kasai2018stochastic}
Hiroyuki Kasai.
\newblock {Stochastic Variance Reduced Multiplicative Update for Nonnegative
  Matrix Factorization}.
\newblock In \emph{2018 IEEE International Conference on Acoustics, Speech and
  Signal Processing (ICASSP)}, pages 6338--6342. IEEE, 2018.

\bibitem[Kim and Klabjan(2019)]{kim2019simple}
Cheolmin Kim and Diego Klabjan.
\newblock {A Simple and Fast Algorithm for L1-norm Kernel PCA}.
\newblock \emph{IEEE Transactions on Pattern Analysis and Machine
  Intelligence}, 2019.

\bibitem[Kim and Klabjan(2020)]{kim2019stochastic}
Cheolmin Kim and Diego Klabjan.
\newblock {Stochastic Variance-Reduced Algorithms for PCA with Arbitrary
  Mini-Batch Sizes}.
\newblock In \emph{International Conference on Artificial Intelligence and
  Statistics}, 2020.

\bibitem[Kim et~al.(2019)Kim, Kim, and Klabjan]{kim2019scale}
Cheolmin Kim, Youngseok Kim, and Diego Klabjan.
\newblock {Scale Invariant Power Iteration}.
\newblock \emph{arXiv preprint arXiv:1905.09882}, 2019.

\bibitem[Kim et~al.(2018)Kim, Carbonetto, Stephens, and Anitescu]{kim2018fast}
Youngseok Kim, Peter Carbonetto, Matthew Stephens, and Mihai Anitescu.
\newblock {A Fast Algorithm for Maximum Likelihood Estimation of Mixture
  Proportions Using Sequential Quadratic Programming}.
\newblock \emph{arXiv preprint arXiv:1806.01412}, 2018.

\bibitem[Lee and Seung(2001)]{lee2001algorithms}
Daniel~D Lee and H~Sebastian Seung.
\newblock {Algorithms for Non-negative Matrix Factorization}.
\newblock In \emph{Advances in Neural Information Processing Systems}, pages
  556--562, 2001.

\bibitem[Lei et~al.(2016)Lei, Zhong, and Dhillon]{lei2016coordinate}
Qi~Lei, Kai Zhong, and Inderjit~S Dhillon.
\newblock {Coordinate-wise Power method}.
\newblock In \emph{Advances in Neural Information Processing Systems}, pages
  2064--2072, 2016.

\bibitem[Lin(2007)]{lin2007projected}
Chih-Jen Lin.
\newblock {Projected Gradient Methods for Non-negative Matrix Factorization}.
\newblock \emph{Neural Computation}, 19\penalty0 (10):\penalty0 2756--2779,
  2007.

\bibitem[Muzzarelli et~al.(2019)Muzzarelli, Weis, Eickhoff, and
  Patil]{muzzarelli2019rank}
Laura Muzzarelli, Susanne Weis, Simon~B Eickhoff, and Kaustubh~R Patil.
\newblock {Rank Selection in Non-negative Matrix Factorization: systematic
  comparison and a new MAD metric}.
\newblock In \emph{2019 International Joint Conference on Neural Networks
  (IJCNN)}, pages 1--8. IEEE, 2019.

\bibitem[Nie et~al.(2016)Nie, Kot{\l}owski, and Warmuth]{nie2016online}
Jiazhong Nie, Wojciech Kot{\l}owski, and Manfred~K Warmuth.
\newblock {Online PCA with Optimal Regret}.
\newblock \emph{Journal of Machine Learning Research}, 17\penalty0
  (1):\penalty0 6022--6070, 2016.

\bibitem[Oja(1982)]{oja1982simplified}
Erkki Oja.
\newblock {Simplified Neuron Model as a Principal Component Analyzer}.
\newblock \emph{Journal of Mathematical Biology}, 15\penalty0 (3):\penalty0
  267--273, 1982.

\bibitem[Oja and Karhunen(1985)]{oja1985stochastic}
Erkki Oja and Juha Karhunen.
\newblock {On Stochastic Approximation of the Eigenvectors and Eigenvalues of
  the Expectation of a Random Matrix}.
\newblock \emph{Journal of Mathematical Analysis and Applications},
  106\penalty0 (1):\penalty0 69--84, 1985.

\bibitem[Serizel et~al.(2016)Serizel, Essid, and Richard]{serizel2016mini}
Romain Serizel, Slim Essid, and Ga{\"e}l Richard.
\newblock {Mini-Batch Stochastic Approaches for Accelerated Multiplicative
  Updates in Nonnegative Matrix Factorisation with Beta-Divergence}.
\newblock In \emph{2016 IEEE 26th International Workshop on Machine Learning
  for Signal Processing (MLSP)}, pages 1--6. IEEE, 2016.

\bibitem[Shamir(2015)]{shamir2015stochastic}
Ohad Shamir.
\newblock {A Stochastic PCA and SVD Algorithm with an Exponential Convergence
  Rate}.
\newblock In \emph{International Conference on Machine Learning}, pages
  144--152, 2015.

\bibitem[Shamir(2016)]{shamir2016convergence}
Ohad Shamir.
\newblock {Convergence of Stochastic Gradient Descent for PCA}.
\newblock In \emph{International Conference on Machine Learning}, pages
  257--265, 2016.

\bibitem[Wang et~al.(2017)Wang, Zhang, and Gu]{wang2017unified}
Lingxiao Wang, Xiao Zhang, and Quanquan Gu.
\newblock {A Unified Variance Reduction-Based Framework for Nonconvex Low-Rank
  Matrix recovery}.
\newblock In \emph{International Conference on Machine Learning}, pages
  3712--3721. PMLR, 2017.

\bibitem[Wang and Zhang(2013)]{wang2013nonnegative}
Yu-Xiong Wang and Yu-Jin Zhang.
\newblock {Nonnegative Matrix Factorization: A Comprehensive Review}.
\newblock \emph{IEEE Transactions on Knowledge and Data Engineering},
  25\penalty0 (6):\penalty0 1336--1353, 2013.

\bibitem[Warmuth and Kuzmin(2008)]{warmuth2008randomized}
Manfred~K Warmuth and Dima Kuzmin.
\newblock {Randomized Online PCA Algorithms with Regret Bounds that are
  Logarithmic in the Dimension}.
\newblock \emph{Journal of Machine Learning Research}, 9\penalty0
  (Oct):\penalty0 2287--2320, 2008.

\bibitem[Xu et~al.(2018)Xu, He, De~Sa, Mitliagkas, and Re]{xu2018accelerated}
Peng Xu, Bryan He, Christopher De~Sa, Ioannis Mitliagkas, and Chris Re.
\newblock {Accelerated Stochastic Power Iteration}.
\newblock In \emph{International Conference on Artificial Intelligence and
  Statistics}, pages 58--67, 2018.

\end{thebibliography}

\newpage
\appendix
%---------------------------------------------%
\section{Description of Data Sets}
\label{sec:data}
%---------------------------------------------%
\begin{table}[ht]
	\caption{Summary of data sets used for KL-NMF} \label{table:1}
	%\vspace{0.1in}
	\centerline{
		\begin{tabular}{|r || r | r | r | r | r |}
		\toprule
		Type & Name & \# of samples & \# of features & \# of nonzeros & Sparsity \\
		\midrule
		Synthetic & Pois1 & 1,000 & 1,000 & 900,000 &
		0.90  \\
		Synthetic & Pois2 & 3,000 & 3,000 & 900,000 &
		0.10  \\
		Synthetic & Pois3 & 9,000 & 9,000 & 900,000 &
		0.01  \\
		\midrule
		Real & NIPS & 1,500 & 12,419 & 280,000 & 
		0.985 \\
		Real & WT & 287 & 19,200 & 5,510,000 &
		0.000 \\
		Real & KOS & 3,430 & 6,906 & 950,000 &
		0.960  \\
		Real & MITF & 361 & 2,429 & 877,000 &
		0.000  \\
		\toprule
		\end{tabular}
		}
\end{table}
%-------------------------------------------------%
The 4 real data sets in the table are retrieved from \url{https://archive.ics.uci.edu/ml/datasets/bag+of+words}, \url{https://www.microsoft.com/en-us/research/project} and \url{https://cbcl.mit.edu/cbcl}. They have already been used in the previous papers such as \cite{hsieh2011fast, kim2019scale}. We preprocess the real data sets by removing few rows and columns having sums less than 20 for NIPS and KOS data sets.

For synthetic data, $V \in \mathbb{R}^{N \times M}$ generated from i.i.d. Poisson random variables, i.e. $V_{ij} \sim \textrm{Poisson}(-\log (1-\rho))$. Here $\rho$ denotes sparsity or proportion of nonzero entries of $V$. This corresponds to the null signal case since in this case KL-NMF is the maximum likelihood estimation problem when $WH = 0$.

%---------------------------------------------%
\section{Proofs}
%---------------------------------------------%
In what follows, we frequently use the fact that for $0 < \eta \leq 1$, $\eta \leq 1/ \max (1, \nu)$ implies 
\begin{align}
    \label{eta-nu-leq-1}
    n \nu \leq 1.
\end{align}
Using $\Delta_0 \leq 1-1/\sqrt{2}$ which follows from \eqref{Delta-0-condition}, we often use
\begin{equation}
    \label{Delta0-upper-bound}
    \frac{\sqrt{\Delta_0}}{1-\Delta_0} \leq 1, \quad \frac{1}{1-\Delta_0} \leq \sqrt{2}.
\end{equation}

\begin{proof}[Proof of Lemma~\ref{lemma:xt-vk-eq}]
%Let $f_{S_t}(x) = \dfrac{1}{|S_t|} \, \sum \limits_{l \in S_t} f_l (x)$. 
From the update rule in Algorithm~\ref{alg:Stochastic-SCI-PI}, we have
\begin{equation}
\label{proof:x-t-recurrence}
\begin{aligned}
    x_{t+1} 
    &= (1-\eta) x_t + \frac{\eta}{\| x_t \|^{p-2}} \left( \nabla f_{S_t}(x_t) - \alpha_t \nabla f_{S_t}(y_0) + \alpha_t \tilde{g} \right)  \\
    &= (1-\eta)  x_t + \frac{\eta}{\| x_t \|^{p-2}} \nabla f(x_t) \\
    & \quad + \frac{\eta}{\| x_t \|^{p-2}} \left[ \nabla f_{S_t}(x_t) - \nabla f(x_t) - \alpha_t \left( \nabla f_{S_t}(y_0) - \nabla f(y_0) \right) \right].
\end{aligned}
\end{equation}
Since $\nabla_i f$ is twice continuously differentiable on an open set containing $\partial \mathcal{B}_{d}$, using the Taylor theorem, we obtain 
%\ys{Taylor theorem}
\begin{align}
\nabla_i f(y_t) = \nabla_i f(x^*) + \nabla \nabla_{i} f(x^*) (y_t - x^*) + \frac{1}{2} \left( y_t -x^* \right)^T H_i(\hat{y}_t^i) \left( y_t - x^* \right)
\label{eq:Taylor-gi}
\end{align}
where 
$\hat{y}_{t}^{i} \in \mathcal{N} (y_t, x^*) \triangleq \{ z \, | \, z = \mu y_t + (1-\mu) x^*, 0 \leq \mu \leq 1 \}$.
Since $f$ is scale invariant with the degree of $p$, by \cite[Proposition 3]{kim2019scale}, we have $c \nabla f(cx) = |c|^{p} \nabla f(x)$, leading to
\begin{align}
\label{proof:grad-xt-vk}
    \frac{\nabla f(x_t)^T z}{\| x_t \|^{p-1}}
    %= \| x_t \|^{p-1} \nabla f(y_t)^T v_k 
    = \nabla f(x^*)^T z + (y_t - x^*)^T \nabla^2 f(x^*) z + \frac{1}{2}  (y_t - x^*)^T \sum_{i=1}^d z_i H_i(\hat{y}_t^i)(y_t - x^*)
\end{align}
for any vector $z \in \mathbb{R}^d$. For $k=1$, using $v_1 = x^*$, we have 
\begin{align*}
    & \nabla f(x^*)^Tv_1 = \nabla f(x^*)^Tx^* = \lambda^*, \\
    & (y_t - x^*)^T \nabla^2 f(x^*) v_1 = (y_t - x^*)^T \nabla^2 f(x^*) x^* = \lambda_1 (y_t^Tx^* - 1),
\end{align*}
which from \eqref{proof:grad-xt-vk} with $z=v_1$ results in
\begin{equation}
\label{proof:grad-xt-v1}
\begin{aligned}
\frac{\nabla f(x_t)^T v_1}{\|x_t\|^{p-1}} & = \lambda^* - \lambda_1 (1-y_t^Tx^*) + \frac{1}{2}  (y_t - x^*)^T \sum_{i=1}^d v_{1i} H_i(\hat{y}_t^i)(y_t - x^*) \\
& = \lambda^* y_t^Tx^* + (\lambda^* - \lambda_1) (1- y_t^Tx^*) + \frac{1}{2}  (y_t - x^*)^T \sum_{i=1}^d v_{1i} H_i(\hat{y}_t^i)(y_t - x^*) \\
& = \lambda^* y_t^Tx^* +  \frac{1}{2}  (y_t - x^*)^T \big[ (\lambda^* - \lambda_1) I + \sum_{i=1}^d v_{1i} H_i(\hat{y}_t^i) \big] (y_t - x^*) \\
& = \lambda^* y_t^Tx^* +  \frac{1}{2}  (y_t - x^*)^T F_1(\hat{y}_t^1,\cdots,\hat{y}_t^d) (y_t - x^*)).
%& = \|x_t\|^{p-2} \left( \lambda^* \|x_t\| - \lambda_1 (\|x_t\|-x_t^Tx^*) + c_{t}^1 (\|x_t\|- x_t^Tx^*) \right) \\
%& = \|x_t\|^{p-2} \left( \lambda^* x_t^Tx^* + (\lambda^* - \lambda_1 + c_{t}^1) (\|x_t\|- x_t^Tx^*) \right).
\end{aligned}
\end{equation}
%Since $x_k^Tx^* \geq 0$, we further have
%\begin{equation}
%\label{proof:xt-norm-xt-v1}
    %\| x_t \| - x_t^Tx^* = x_t^Tx^* \left( \frac{\|x_t\|}{x_t^Tx^*} - 1 \right) \leq
    %x_t^Tx^* \sqrt{\frac{\sum_{k=2}^d (x_t^Tv_k)^2}{(x_t^Tx^*)^2}} \leq x_t^Tx^* %\sqrt{\frac{1-(y_t^Tx^*)^2}{(y_t^Tx^*)^2}}.
%\end{equation}
For $2 \leq k \leq d$, from \eqref{proof:grad-xt-vk} with $z=v_k$, $(x^*)^Tv_k = v_1^T v_k = 0$ and $\nabla f(x^*)^Tv_k = \lambda^* v_1^Tv_k = 0$, we have
\begin{equation}
\label{proof:grad-f-v-k}
\begin{aligned}
\frac{\nabla f(x_t)^T v_k}{\|x_t\|^{p-1}} & = \lambda_k y_t^Tx^* +  \frac{1}{2}  (y_t - x^*)^T F_k(\hat{y}_t^1,\cdots,\hat{y}_t^d) (y_t - x^*).
\end{aligned}
\end{equation}
Since $\nabla f_l$ is scale invariant with the degree of $p-1$ for each $l \in [n]$, we have 
\begin{equation*}
    \nabla f_l(x_t) = \|x_t\|^{p-1} \nabla f_l(y_t), \quad \alpha_t \nabla f_l(y_0) = \|x_t\|^{p-1} (y_t^T y_0)^{p-1} \nabla f_l(y_0),
\end{equation*}
which leads to
\begin{align*}
    \frac{1}{\|x_t \|^{p-1}} \left( \nabla f_{S_t}(x_t) - \nabla f (x_t) - \alpha_t \left( \nabla f_{S_t}(y_0) - \nabla f(y_0) \right) \right)
    = \nabla g_{S_t} (y_t) - \nabla g_{S_t} \big( (y_t^Ty_0) y_0 \big).
\end{align*}
%Therefore, we have
%\begin{equation} 
%    \frac{1}{\| x_t \|^{p-2}} \sum_{k=1}^d \big[ \nabla_{S_t} f(x_t) - \nabla f(x_t) - \alpha_t \nabla_{S_t} f(x_0) - \nabla f(x_0)) \big]^T v_k 
%    =    - (x_t - (x_t^Ty_0)y_0) \left( \sum_{k=1}^d G_t^T v_k \right).
%\end{equation}
Using the Taylor approximation of $\nabla_k g_{S_t}$ around $(y_t^Ty_0)y_0$, we have
\begin{align*}
    \nabla_k g_{S_t}(y_t) - \nabla_k g_{S_t}((y_t^Ty_0)y_0) = \nabla \nabla_k g_{S_t} \big( \bar{y}_t^k \big)^T \big( y_t-(y_t^Ty_0) y_0 \big)
\end{align*}
where $\bar{y}_{t}^{k} \in \mathcal{N} (y_t, (y_t^Ty_0) y_0)$. This leads to
\begin{align}
\label{proof:grad-xt-alpha-grad-x0}
    \frac{1}{\|x_t \|^{p-2}} \left( \nabla f_{S_t}(x_t) - \nabla f (x_t) - \alpha_t \left( \nabla f_{S_t}(y_0) - \nabla f(y_0) \right) \right)
    =  G_{S_t}(\bar{y}_t^1,\cdots,\bar{y}_t^d) \big( x_t-(x_t^Ty_0)y_0  \big).
\end{align}
Using \eqref{proof:x-t-recurrence}, \eqref{proof:grad-xt-v1}, \eqref{proof:grad-f-v-k} and \eqref{proof:grad-xt-alpha-grad-x0}, we have
\begin{equation}
\label{proof:x-t-1-v-1}
\begin{aligned}
    x_{t+1}^T v_k &= \left(1-\eta+ \eta (\lambda_k + (\lambda^* - \lambda_1) \mathds{1}_{k = 1}) \right) x_t^Tv_k  +  \frac{1}{2} \eta \| x_t \| (y_t - x^*)^T F_k(\hat{y}_t^1,\cdots,\hat{y}_t^d) (y_t - x^*) \\
    & \qquad + \eta \left( G_{S_t}(\bar{y}_t^1,\cdots,\bar{y}_t^d) \big( x_t-(x_t^Ty_0)y_0  \big) \right)^T v_k.
\end{aligned}
\end{equation}
\end{proof}

\begin{proof}[Proof of Lemma~\ref{lemma:y0-St-necessary-condition}]
We prove by induction. Suppose that we have $\Delta_s \leq \Delta_0$ for $s \leq t < m$. Since $\Delta_0 \leq 1-1/\sqrt{2}$, this implies that $y_t^Tx^* \geq 1/\sqrt{2}$ and $y_0^Tx^* \geq 1/\sqrt{2}$. Therefore, we have
\begin{align*}
    y_t^Ty_0 &= \left[ (y_t^Tx^*)x^* + y_t - (y_t^Tx^*)x^* \right]^T \left[ (y_0^Tx^*)x^* + y_0 - (y_0^Tx^*)x^* \right] \\
    &= (y_t^Tx^*)(y_0^Tx^*) + (y_t - (y_t^Tx^*)x^*)^T (y_0 - (y_0^Tx^*)x^*) \\
    &\geq (y_t^Tx^*)(y_0^Tx^*) - \| y_t - (y_t^Tx^*)x^* \| \| y_0 - (y_0^Tx^*)x^* \| \\
    &\geq (y_t^Tx^*)(y_0^Tx^*) - \sqrt{1-(y_t^Tx^*)^2} \sqrt{1-(y_0^Tx^*)^2} \\
    &\geq 0,
\end{align*}
which leads to
\begin{align*}
    \| x_t - (x_t^Ty_0)y_0 \|^2 = \| x_t \|^2 (1 - (y_k^T y_0)^2) &\leq 2 \| x_t \|^2 (1- y_t^T y_0) = \|x_t \|^2 \| y_t - y_0 \|^2.
\end{align*}
By the triangular inequality, $(a+b)^2 \leq 2(a^2+b^2)$ and $\Delta_t \leq \Delta_0$, we have
\begin{align*}
    \| y_t - y_0 \|^2 \leq 2( \|y_t - x^* \|^2 + \|y_0 - x^* \|^2) \leq 4  \|y_0 - x^* \|^2.
\end{align*}
From $y_0^Tx^* \geq 0$, we further obtain
\begin{align}
    \| x_t - (x_t^Ty_0)y_0 \|^2 &\leq 4 \| x_t \|^2 \|y_0 - x^* \|^2 = 8 \|x_t \|^2 (1- y_0^Tx^*) \label{eq:yt-y0-bound-1} \\
    &\leq 8 \|x_t\|^2 (1-(y_0^Tx^*)^2) = 8 \|x_t\|^2 \sum_{k=2}^d (y_0^Tv_k)^2.
    \label{eq:yt-y0-bound}
\end{align}
Using Lemma~\ref{lemma:xt-vk-eq}, the definitions of $M$ and $L$, \eqref{eq:yt-y0-bound-1} and that $\Delta_t \leq \Delta_0$, we have
\begin{equation}
    \label{eq:xt-v1-square-bound-L}
    \begin{aligned}
    x_{t+1}^T v_1 & \geq \left(1-\eta+ \eta \lambda^* \right) x_t^Tv_1  -  \frac{1}{2} \eta M \| x_t \| \|y_t - x^*\|^2 - \eta \sqrt{L} \| x_t-(x_t^Ty_0)y_0 \| \\
    & \geq \left(1-\eta+ \eta \lambda^* \right) x_t^Tv_1  - \eta M (1-y_t^Tx^*) \| x_t \| -  \eta \sqrt{8L(1-y_0^Tx^*)} \|x_t \| \\
    & \geq \left[ 1-\eta+ \eta \left( \lambda^* - \frac{M \Delta_0}{1-\Delta_0} - \frac{\sqrt{8L \Delta_0}}{1-\Delta_0} \right) \right] y_0^Tx^* \|x_t \|.
    \end{aligned}
\end{equation}
By \eqref{Delta0-upper-bound}, \eqref{Delta-0-condition} and that $L \leq L_0$, we have
\begin{align*}
    \lambda^* - \frac{M \Delta_0}{1-\Delta_0} - \frac{\sqrt{8L \Delta_0}}{1-\Delta_0}
    \geq \lambda^* - \big( M + 4 \sqrt{L} \big)  \sqrt{\Delta_0} 
    \geq \lambda^* - \frac{(\lambda^* - \bar{\lambda}) \big( M + 4 \sqrt{L} \big)}{2M + 4 \sqrt{L_0}} 
    \geq 0.
\end{align*}
This leads to $x_{t+1}^Tv_1 \geq 0$.

Now, we prove that $\Delta_{t+1} \leq \Delta_0$. Since $\{v_1,\ldots,v_d \}$ forms an orthogonal basis, we have $\| x_t\|^2 = \sum_{k=1}^d (x_t^Tv_k)^2$.
Since
\begin{align}
    &\sum_{k=2}^d (1-\eta+\eta \lambda_k)^2 (x_t^Tv_k)^2 \leq (1-\eta+\eta \bar{\lambda})^2 \sum_{k=2}^d (x_t^Tv_k)^2 \label{eq:xt-vk-first-bound}
    \\
    &\sum_{k=1}^d \left(1-\eta+ \eta (\lambda_k + (\lambda^* - \lambda_1) \mathds{1}_{k = 1}) \right)^2 (x_t^Tv_k)^2 \leq (1-\eta+\eta \lambda^*)^2 \| x_t \|^2 \label{eq:xt-norm-first-bound}
\end{align}
\vspace{-2mm}
\begin{equation}
\label{eq:xt-vk-second-bound}
\begin{aligned}
    \sum_{k=2}^d \left[(y_t - x^*)^T F_k(\hat{y}_t^1,\cdots,\hat{y}_t^d) (y_t - x^*) \right]^2 
    &\leq \sum_{k=1}^d \left[(y_t - x^*)^T F_k(\hat{y}_t^1,\cdots,\hat{y}_t^d) (y_t - x^*) \right]^2 \\
    & \leq M^2 \|y_t - x^*\|^4
\end{aligned}
\end{equation}
\vspace{-2mm}
\begin{equation}
\label{eq:xt-vk-third-bound}
\begin{aligned}
    \sum_{k=2}^d \left[ v_k^T G_{S_t}(\bar{y}_t^1,\cdots,\bar{y}_t^d) \big( x_t-(x_t^Ty_0)y_0  \big) \right]^2 & \leq \sum_{k=1}^d \left[ v_k^T G_{S_t}(\bar{y}_t^1,\cdots,\bar{y}_t^d) \big( x_t-(x_t^Ty_0)y_0 \big) \right]^2 \\
    & \leq L \| x_t - (x_t^Ty_0) y_0 \|^2
\end{aligned}
\end{equation}
where \eqref{eq:xt-vk-third-bound} follows from $\| \sum_{k=1}^d v_kv_k^T \| = 1$. By Lemma~\ref{lemma:xt-vk-eq} and the Cauchy-Schwarz inequality, we have
\begin{align}
    & \sum_{k=2}^d (x_{t+1}^T v_k)^2 
    \leq \big[ (1-\eta+\eta \bar{\lambda}) \sqrt{\sum_{k=2}^d (x_t^Tv_k)^2} + \frac{1}{2} \eta M \|x_t\| \|y_t-x^*\|^2 + \eta \sqrt{L} \|x_t - (x_t^Ty_0)y_0 \| \big]^2 \label{eq:xt-vk-square-sum-relation} \\
    & \| x_{t+1} \|^2 = \sum_{k=1}^d (x_{t+1}^Tv_k)^2 
    \leq \big[ 1-\eta+\eta \lambda^* + \frac{1}{2} \eta M \| y_0 - x^* \|^2 + \eta \sqrt{L} \|y_t - (y_t^Ty_0)y_0 \|  \big]^2 \| x_t \|^2. \label{eq:xt-norm-square-relation}
\end{align}
First, we consider the case when \eqref{L-condition} holds. From $\Delta_t \leq \Delta_0 \leq 1$, we have $0 \leq y_t^Tx^* \leq 1$ and $\sum_{k=2}^d (y_t^Tv_k)^2 = 1 - (y_t^Tx^*)^2 \leq 1 - (y_0^Tx^*)^2 = \sum_{k=2}^d (y_0^Tv_k)^2$, resulting in
\begin{equation}
\label{eq:yt-opt-gap-bound}
\begin{aligned}
\| y_t - x^* \|^2 \leq 2 \sqrt{1-y_t^Tx^*} \sqrt{1-(y_t^Tx^*)^2} \leq 2 {\sqrt{\Delta_t}} \sqrt{\sum_{k=2}^d (y_0^Tv_k)^2}
\leq 2 \sqrt{\Delta_0} \sqrt{\sum_{k=2}^d (y_0^Tv_k)^2}.
\end{aligned}
\end{equation}
Plugging \eqref{eq:yt-y0-bound} and \eqref{eq:yt-opt-gap-bound} into \eqref{eq:xt-vk-square-sum-relation}, we have
\begin{align}
    \label{eq:xt-vk-square-sum-bound-L}
    \sum_{k=2}^d (x_{t+1}^T v_k)^2 
    &\leq 
    \left[ 1-\eta+\eta \left( \bar{\lambda} + M \sqrt{\Delta_0} + 2 \sqrt{2L} \right)  \right]^2 \| x_t \|^2 \sum_{k=2}^d (y_0^Tv_k)^2.
\end{align}
Combining \eqref{eq:xt-v1-square-bound-L} and \eqref{eq:xt-vk-square-sum-bound-L}, we have
\begin{align}
    \label{eq:lemma-3-2-final-1}
    \frac{\sum_{k=2}^d (x_{t+1}^T v_k)^2 }{(x_{t+1}^T v_1)^2}
    \leq
    \bigg[
    \frac{1-\eta+\eta \big( \bar{\lambda} + M \sqrt{\Delta_0} + 2 \sqrt{2L} \big) }
    {1-\eta+ \eta \left( \lambda^* - M \Delta_0/(1-\Delta_0) - 2\sqrt{2L \Delta_0}/(1-\Delta_0) \right)}
    \bigg]^{2}
    \frac{\sum_{k=2}^d (y_{0}^T v_k)^2 }{(y_{0}^T v_1)^2}.
\end{align}
Using \eqref{Delta0-upper-bound} and \eqref{L-condition}, we have 
\begin{align*}
    \lambda^* - \frac{M\Delta_0}{1-\Delta_0} -  \frac{2 \sqrt{2L \Delta_0}}{1-\Delta_0}
    - \big( \bar{\lambda} + M \sqrt{\Delta_0} + 2\sqrt{2L} \big)
    \geq
    (\lambda^* - \bar{\lambda})
    - 2M \sqrt{\Delta_0}
    - 4\sqrt{2L}
    \geq 0.
\end{align*}
Therefore, from \eqref{eq:lemma-3-2-final-1}, we finally have 
\begin{equation*}
    \frac{1 - (y_{t+1}^Tx^*)^2}{(y_{t+1}^Tx^*)^2} = \frac{\sum_{k=2}^d (y_{t+1}^T v_k)^2 }{(y_{t+1}^T v_1)^2} = \frac{\sum_{k=2}^d (x_{t+1}^T v_k)^2 }{(x_{t+1}^T v_1)^2}
    \leq \frac{\sum_{k=2}^d (y_{0}^T v_k)^2 }{(y_{0}^T v_1)^2} = \frac{1 - (y_{0}^Tx^*)^2}{(y_{0}^Tx^*)^2},
\end{equation*}
which leads to $\Delta_{t+1} = 1-y_{t+1}^Tx^* \leq 1 - y_0^Tx^* = \Delta_0$.

Next, we derive $\Delta_{t+1} \leq \Delta_0$ from \eqref{eta-condition}. From \eqref{eq:yt-y0-bound-1} and \eqref{eq:xt-norm-square-relation}, we have
\begin{align*}
    \| x_{t+1} \|^2
    &\leq \Big[ 1-\eta+\eta \Big(
    \lambda^* + \frac{1}{2} M \| y_0 - x^* \|^2 + 2 \sqrt{L} \| y_0 - x^* \| \Big) \Big]^2 \| x_t \|^2.
\end{align*}
Using induction, this leads to
\begin{align}
    \| x_{t+1} \|^2
    &\leq \Big[ 1-\eta+\eta \Big( \lambda^* + \frac{1}{2} M \| y_0 - x^* \|^2 + 2 \sqrt{L} \| y_0 - x^* \| \Big)  \Big]^{2(t+1)} \| x_0 \|^2. \label{eq:xt-norm-square-relation-bound}
\end{align}
On the other hand, from \eqref{proof:x-t-recurrence}, \eqref{proof:grad-xt-alpha-grad-x0}, \eqref{eq:yt-y0-bound-1} and the definition of $L$, we have
\begin{align*}
    x_{t+1}^Ty_0 &= (1-\eta) x_{t}^T y_0 + \frac{\eta \nabla f(x_{t})^Ty_0}{\| x_{t} \|^{p-2}} + \eta y_0^T G_{S_{t}} (\bar{y}_{t}^1, \cdots, \bar{y}_{t}^d) (x_{t} - (x_{t}^Ty_0) y_0) \\
    &\geq (1-\eta) x_{t}^T y_0 + \frac{\eta \nabla f(x_{t})^Ty_0}{\| x_{t} \|^{p-2}}
    - 2\eta \sqrt{L} \| y_0 - x^*\| \| x_{t} \|.
\end{align*}
Using $z=y_0$ in \eqref{proof:grad-xt-vk} and using $\nabla f(x^*) = \lambda^* x^*$ and the definition of $M$, we have
\begin{align*}
    \frac{\nabla f(x_{t})^Ty_0}{\| x_t \|^{p-1}} 
    &= \nabla f(x^*)^T y_0 + (y_t-x^*)^T \nabla^2 f(x^*) y_0 + \frac{1}{2} (y_t-x^*)^T \sum_{i=1}^d y_{0i} H_i(\hat{y}_t^i) (y_t -x^*)  \\
    &= \lambda^* y_{t}^Ty_0 + (y_t-x^*)^T \left( \nabla^2 f(x^*) - \lambda^* I \right) y_0 - \frac{1}{2} M \| y_t - x^* \|^2 \\
    &\geq \lambda^* y_{t}^Ty_0 - (\lambda^* + \sigma) \| y_t - x^* \| - \frac{1}{2} M \| y_t - x^* \|^2 \\
    &\geq \lambda^* y_t^Ty_0 - (\lambda^* + \sigma) \| y_0 - x^* \| - \frac{1}{2} M \| y_0 - x^* \|^2
\end{align*}
where the last inequality follows from
\begin{equation*}
    \| y_t - x^* \|^2 = 2(1-y_t^Tx^*) \le 2(1-y_0^Tx^*) = \| y_0 - x^* \|^2.
\end{equation*}
This results in
\begin{align*}
        x_{t+1}^Ty_0 &\geq (1-\eta+\eta \lambda^*) x_{t}^Ty_0 - \eta \big( \lambda^* + \sigma + \frac{1}{2} M \| y_0 - x^* \| + 2\sqrt{L} \big) \| y_0 - x^*\| \| x_{t} \| \\
        &= (1-\eta+\eta \lambda^*) x_{t}^Ty_0 - \eta \theta_1 \sqrt{2 \Delta_0} \| x_t \|.
\end{align*}
Combining with \eqref{eq:xt-norm-square-relation-bound}, we obtain
\begin{align*}
    x_{t+1}^Ty_0 
    &\geq (1-\eta+\eta \lambda^*) x_{t}^Ty_0 \\
    & \quad - \eta \theta_1 \sqrt{2 \Delta_0} \Big[ 1-\eta+\eta \Big( \lambda^* + \frac{1}{2} M \| y_0 - x^* \|^2 + 2 \sqrt{L} \| y_0 - x^* \| \Big)  \Big]^{t} \| x_0 \| \\
    &\geq (1-\eta+\eta \lambda^*) x_{t}^Ty_0  - \eta \theta_1 \sqrt{2 \Delta_0} \Big[ 1-\eta+\eta \lambda^* + \eta \theta_1 \sqrt{2\Delta_0}  \Big]^{t} \| x_0 \|.
\end{align*}
By recursion, we further have
\begin{align}
    x_{t+1}^Ty_0
    &\geq (1-\eta+\eta \lambda^*)^{t+1} \| x_0 \| \\
    &\quad - \eta \theta_1 \sqrt{2 \Delta_0} \sum_{i=1}^{t+1} \big( 1-\eta+\eta \lambda^* \big)^{i-1} \Big[ 1-\eta+\eta \lambda^* + \eta \theta_1 \sqrt{2 \Delta_0} \Big]^{t+1-i} \| x_0 \| \nonumber\\
    &= \left[ 2 (1-\eta+\eta \lambda^*)^{t+1} -  \left( 1-\eta+\eta \lambda^* + \eta \theta_1 \sqrt{2 \Delta_0} \right)^{t+1} \right] \| x_0 \|. \label{eq:xt-y0-lb}
\end{align}
Also, by the definition of $\nu_1$ and requirement \eqref{eta-condition} that $\eta \nu_1 \leq 1$ which yields
\begin{equation*}
    \bar{x} = \frac{\eta \theta_1 \sqrt{2\Delta_0}}{1-\eta+\eta \lambda^*} \leq \frac{1}{2m},
\end{equation*}
it is easy to establish that 
\begin{equation}
    \label{proof:poly-exp-linear}
    \left( 1+x \right)^t \leq \exp \left(xt \right) \leq 2xt+1
\end{equation}
for any $0 \leq x \leq 1/2t$. Since $t < m$, by considering $x= \bar{x}$, we obtain from \eqref{eq:xt-y0-lb} inequality
\begin{equation*}
    x_{t+1}^Ty_0 \geq (1-\eta+\eta \lambda^*)^{t+1}(1-2m\bar{x}) \geq 0.
\end{equation*}
Since $\| x_t - (x_t^Ty_0)y_0 \|^2 = \| x_t \|^2 - (x_t^Ty_0)^2$, using \eqref{eq:xt-norm-square-relation-bound}, \eqref{eq:xt-y0-lb} and elementary algebraic manipulations, we have
\begin{equation*}
\begin{aligned}
    \| x_t - (x_t^Ty_0)y_0 \|^2 
    &\leq 4 (1-\eta+\eta \lambda^*)^{2t} \Big[ \Big( 1+\frac{\eta \theta_1 \sqrt{2\Delta_0}}{1-\eta+\eta \lambda^*} \Big)^{t} - 1 \Big] \| x_0 \|^2.
\end{aligned}
\end{equation*}
By \eqref{eta-condition}, \eqref{def:nu-1} and \eqref{eta-nu-leq-1}, we have $\eta (1-\lambda^* + 2 \theta_1 m \sqrt{2 \Delta_0}) \leq 1$ or 
\begin{equation*}
    \frac{\eta \theta_1 m \sqrt{2\Delta_0}}{1-\eta+\eta \lambda^*} \leq \frac{1}{2}.
\end{equation*}
Since
\begin{equation*}
    \frac{\eta \theta_1 t \sqrt{2 \Delta_0}}{1-\eta+\eta \lambda^*} 
    \leq
    \frac{\eta \theta_1 m \sqrt{2 \Delta_0}}{1-\eta+\eta \lambda^*} 
    \leq \frac{1}{2} < 1,
\end{equation*}
using \eqref{proof:poly-exp-linear}, we obtain
\begin{equation}
\label{eq:xt-norm-diff-xt-y0}
    \| x_t - (x_t^Ty_0)y_0 \|^2 \leq 8 \eta \theta_1 (1-\eta+\eta \lambda^*)^{2t-1}  \sqrt{2\Delta_0} t \| x_0 \|^2.
\end{equation}
Plugging \eqref{eq:yt-opt-gap-bound} and \eqref{eq:xt-norm-diff-xt-y0} into the square root of \eqref{eq:xt-vk-square-sum-relation} and then apply recursion, we have
\begin{equation}
\label{eq:xt-vk-square-root-bound}
\begin{aligned}
    \sqrt{\sum_{k=2}^d (x_{t+1}^Tv_k)^2} 
    &\leq \big[ 1-\eta+\eta (\bar{\lambda}+M\sqrt{\Delta_0}) \big] \sqrt{\sum_{k=2}^d (x_{t}^Tv_k)^2} \\
    & \quad + \eta \sqrt{\frac{8 \eta L \theta_1 \sqrt{2\Delta_0}}{1-\eta+\eta \lambda^*}} (1-\eta+\eta \lambda^*)^t t \| x_0 \| \\
    &\leq \big[ 1-\eta+\eta (\bar{\lambda}+M\sqrt{\Delta_0}) \big]^{t+1} \sqrt{\sum_{k=2}^d (x_{0}^Tv_k)^2} \\
    &\quad + \eta \sqrt{\frac{8 \eta L \theta_1 \sqrt{2\Delta_0}}{1-\eta+\eta \lambda^*}} \sum_{i=1}^t i \, (1-\eta+\eta \lambda^*)^i \big[ 1-\eta+\eta (\bar{\lambda}+M\sqrt{\Delta_0}) \big]^{t-i} \| x_0 \|.
\end{aligned}
\end{equation}
For a positive integer $t$ and a non-negative real number $r \geq 0$ such that $rt \leq 1$, we have
\begin{align*}
    (1+r)^t - 1 &= r \left( (1+r)^{t-1} + (1+r)^{t-2} + \cdots + 1 \right) \geq rt
    %(1+r)^t - 1 &\leq \exp(rt) - 1 \leq 2rt,
\end{align*}
and \eqref{proof:poly-exp-linear} with $x=r$, which results in
\begin{equation}
\label{eq:sum-i-ri-bound}
\begin{aligned}
    \sum_{i=1}^t (1+r)^i i &= \frac{1+r}{r^2} \left( t (1+r)^{t+1} - (t+1) (1+r)^t + 1 \right) \\
    &\leq \frac{1+r}{r^2} \left( t (1+r)^{t+1} - t (1+r)^t - rt \right) \\
    &= \frac{(1+r)t}{r} \left( (1+r)^{t} - 1 \right) \\
    &\leq 2(1+r) t^2.
\end{aligned}
\end{equation}
By \eqref{eta-condition}, \eqref{def:nu-2} and \eqref{eta-nu-leq-1}, we have
\begin{equation*}
    \frac{1-\eta+\eta \lambda^*}{1-\eta+\eta (\bar{\lambda} + M \sqrt{\Delta_0})} - 1 \leq \frac{1}{m}.
\end{equation*}
Also, by \eqref{Delta-0-condition}, we have $\lambda^* - \bar{\lambda} - M \sqrt{\Delta_0} \geq 0$, leading to
\begin{equation*}
    \frac{1-\eta+\eta \lambda^*}{1-\eta+\eta (\bar{\lambda} + M \sqrt{\Delta_0})} - 1 = \frac{\eta (\lambda^* - \bar{\lambda} - M \sqrt{\Delta_0})}{1-\eta+\eta (\bar{\lambda} + M \sqrt{\Delta_0})} \geq 0.
\end{equation*}
Therefore, using \eqref{eq:sum-i-ri-bound}, we have
\begin{equation}
\label{eq:xt-vk-square-root-bound-second}
\begin{aligned}
    & \sum_{i=1}^t i \, (1-\eta+\eta \lambda^*)^i \big[ 1-\eta+\eta (\bar{\lambda}+M\sqrt{\Delta_0}) \big]^{t-i} \\
    &= \left[ 1-\eta+\eta (\bar{\lambda}+M\sqrt{\Delta_0}) \right]^{t} \sum_{i=1}^t i \, \left[ \frac{1-\eta+\eta \lambda^*}{1-\eta+\eta(\bar{\lambda}+M \sqrt{\Delta_0})} \right]^{i} \\
    & \leq 2 (1-\eta+\eta \lambda^*) t^2 \big[ 1-\eta+\eta (\bar{\lambda}+M\sqrt{\Delta_0}) \big]^{t-1}.
\end{aligned}
\end{equation}
Plugging \eqref{eq:xt-vk-square-root-bound-second} into \eqref{eq:xt-vk-square-root-bound}, we obtain
\begin{equation}
\label{eq:xt-vk-square-root-bound-final}
\begin{aligned}
    \sqrt{\sum_{k=2}^d (x_{t+1}^Tv_k)^2}
    &\leq 
    \big[ 1-\eta+\eta (\bar{\lambda}+M\sqrt{\Delta_0}) \big]^{t+1} \sqrt{\sum_{k=2}^d (x_{0}^Tv_k)^2} \\
    &\quad + 2 \eta \sqrt{8 (1-\eta+\eta \lambda^*) \eta L \theta_1 \sqrt{2\Delta_0}} \, t^2
      \big[ 1-\eta+\eta (\bar{\lambda}+M\sqrt{\Delta_0}) \big]^{t-1} \| x_0 \|.
\end{aligned}
\end{equation}
On the other hand, from \eqref{eq:xt-v1-square-bound-L} and
\begin{align*}
    (1-y_t^Tv_1) \| x_t \| &= \frac{1-y_t^Tv_1}{y_t^Tv_1} x_t^Tv_1 \leq \frac{1-y_0^Tv_1}{y_0^Tv_1} x_t^Tv_1 = \frac{\Delta_0}{1-\Delta_0} x_t^Tv_1, \\
    \sqrt{1-y_0^Tv_1} \| x_t \| &= \frac{\sqrt{1-y_0^Tv_1}}{y_t^Tv_1} x_t^Tv_1 \leq \frac{\sqrt{1-y_0^Tv_1}}{y_0^Tv_1} x_t^Tv_1 = \frac{\sqrt{1-\Delta_0}}{1-\Delta_0} x_t^Tv_1,
\end{align*}
we have
\begin{equation}
    \label{eq:xt-v1-square-bound-L-epoch}
    x_{t+1}^Tv_1 \geq \left[ 1 - \eta + \eta \left( \lambda^* - \frac{M \Delta_0}{1- \Delta_0} - \frac{2\sqrt{2L \Delta_0}}{1-\Delta_0} \right) \right]^{t+1} x_0^Tv_1.
\end{equation}
Combining \eqref{eq:xt-vk-square-root-bound-final} and \eqref{eq:xt-v1-square-bound-L-epoch}, we have
\begin{equation}
\label{eq:lemma-3-2-bound}
\begin{aligned}
    \frac{\sqrt{\sum_{k=2}^d (x_{t+1}^Tv_k)^2}}{x_{t+1}^Tv_1}
    & \leq 
    \Bigg[ \frac{1-\eta+\eta (\bar{\lambda}+M\sqrt{\Delta_0})}{1 - \eta + \eta \left[ \lambda^* - ({M \Delta_0}+{2\sqrt{2L \Delta_0}})/(1-\Delta_0) \right]} \Bigg]^{t+1} \frac{\sqrt{\sum_{k=2}^d (x_{0}^Tv_k)^2}}{x_{0}^Tv_1} \\
    & \quad + \frac{2 \eta t^2 \sqrt{8 (1-\eta+\eta \lambda^*) \eta L \theta_1 \sqrt{2\Delta_0}} \,
      \big[ 1-\eta+\eta (\bar{\lambda}+M\sqrt{\Delta_0}) \big]^{t-1}}{\left( 1 - \eta + \eta \left[ \lambda^* - ({M \Delta_0}+{2\sqrt{2L \Delta_0}})/(1-\Delta_0) \right] \right)^{t+1} y_0^Tv_1}.
\end{aligned}
\end{equation}
Since $0 < \eta \leq 1$ and $\bar{\lambda} < \lambda^*$, we have
\begin{equation}
\label{eq:lemma-3-2-final-2}
    \frac{\bar{\lambda}}{\lambda^*} \leq \frac{1-\eta+\eta \bar{\lambda}}{1-\eta+\eta \lambda^*}
    \leq \frac{1-\eta+\eta (\bar{\lambda}+M\sqrt{\Delta_0})}{1 - \eta + \eta \lambda^*}.
\end{equation}
Let
\begin{equation}
\label{def:gamma}
    \gamma = \frac{\lambda^*-\bar{\lambda}-M\sqrt{\Delta_0}-(M\Delta_0 + 2\sqrt{2L\Delta_0})/(1-\Delta_0)}{1 - \eta + \eta \left[ \lambda^* - ({M \Delta_0}+{2\sqrt{2L \Delta_0}})/(1-\Delta_0) \right]}.
\end{equation}
By \eqref{Delta0-upper-bound} and $\theta_2 \geq 0$ due to \eqref{Delta-0-condition}, we have
\begin{equation}
\label{eq:lemma-3-2-final-2-2}
\begin{aligned}
    \frac{1}{1-\eta+\eta (\bar{\lambda}+M\sqrt{\Delta_0})}
    &= \frac{\gamma}{1-\eta \gamma} \bigg[ \frac{1}{\lambda^*-\bar{\lambda}-M\sqrt{\Delta_0}-(M\Delta_0 + 2\sqrt{2L\Delta_0})/(1-\Delta_0)} \bigg] \\
    &\leq \frac{\gamma}{\theta_2 (1-\eta \gamma)}.
\end{aligned}
\end{equation}
Using \eqref{eq:lemma-3-2-final-2}, \eqref{eq:lemma-3-2-final-2-2} and that $y_0^Tv_1 \geq 1/\sqrt{2}$, we have
\begin{align*}
    \frac{2 \eta t^2 \sqrt{8 (1-\eta+\eta \lambda^*) \eta L \theta_1 \sqrt{2\Delta_0}}}{y_0^Tv_1 (1-\eta+\eta(\bar{\lambda}+M\sqrt{\Delta_0}))^2}
    &\leq 8 \sqrt{2} \sqrt{\frac{\lambda^*}{\bar{\lambda}}} \sqrt{\frac{\eta L \theta_1 \sqrt{\Delta_0}}{1-\eta+\eta (\bar{\lambda}+M\sqrt{\Delta_0})}}
    \frac{\eta \gamma t^2}{\theta_2(1-\eta \gamma)}.
\end{align*}
By \eqref{eta-condition}, \eqref{def:nu-3} and \eqref{eta-nu-leq-1}, we have
\begin{align*}
    \eta \left( \frac{128 L \theta_1 \lambda^* m^2}{\theta_2^2 \bar{\lambda} \Delta_0 \sqrt{\Delta_0}} + 1 - \big( \bar{\lambda} + M \sqrt{\Delta_0} \big) \right) \leq 1
\end{align*}
or
\begin{align*}
    \frac{\eta L \theta_1 \sqrt{\Delta_0}}{1-\eta+\eta \big( \bar{\lambda} + M \sqrt{\Delta_0}\big)}
    \leq \frac{\theta_2^2 \bar{\lambda} \Delta_0^2}{128 \lambda^* m^2},
\end{align*}
which results in
\begin{align}
\label{eq:lemma-3-2-final-3}
    \frac{2 \eta t^2 \sqrt{8 (1-\eta+\eta \lambda^*) \eta L \theta_1 \sqrt{2\Delta_0}}}{y_0^Tv_1 (1-\eta+\eta(\bar{\lambda}+M\sqrt{\Delta_0}))^2}
    &\leq \frac{\eta \gamma t^2 \Delta_0}{(1-\eta \gamma)m} \leq \frac{\eta \gamma t^2}{(1-\eta \gamma)m} \frac{\sum_{k=2}^d (x_0^Tv_k)^2}{(x_0^Tv_1)^2}.
\end{align}
The last inequality follows from
\begin{equation*}
    \Delta_0 = 1-y_0^Tx^* \leq 1-(y_0^Tx^*)^2 \leq \frac{\sum_{k=2}^d (y_0^Tv_k)^2}{(y_0^Tv_1)^2} = \frac{\sum_{k=2}^d (x_0^Tv_k)^2}{(x_0^Tv_1)^2}.
\end{equation*}
Plugging \eqref{def:gamma} and \eqref{eq:lemma-3-2-final-3} into \eqref{eq:lemma-3-2-bound}, we have
\begin{align*}
    \frac{\sqrt{\sum_{k=2}^d (x_{t+1}^Tv_k)^2}}{x_{t+1}^Tv_1} 
    &\leq 
    (1-\eta \gamma)^{t+1} \bigg[ 1 
    + \frac{\eta \gamma t^2}{(1-\eta \gamma)m} \bigg] \frac{\sqrt{\sum_{k=2}^d (x_0^Tv_k)^2}}{x_0^Tv_1}.
\end{align*}
Using $1+nx \leq (1+x)^n$ for $x \geq 0$ and the fact that $\gamma \geq 0$ by \eqref{Delta-0-condition}, we have
\begin{align*}
    (1-\eta \gamma)^{t+1} \bigg[ 1 
    + \frac{\eta \gamma t^2}{(1-\eta \gamma)m} \bigg]
    &= 1 - \bigg[ \bigg(1+\frac{\eta \gamma}{1-\eta \gamma}\bigg)^{t+1} - 1 - \frac{\eta \gamma t^2}{(1-\eta \gamma)m} \bigg] (1-\eta \gamma)^{t+1} \\
    &\leq 1 - \bigg( t+1  - \frac{t^2}{m} \bigg) \eta \gamma (1-\eta \gamma)^{t},
\end{align*}
which yields
\begin{align*}
    \frac{\sqrt{\sum_{k=2}^d (x_{t+1}^Tv_k)^2}}{x_{t+1}^Tv_1} \leq   \frac{\sqrt{\sum_{k=2}^d (x_0^Tv_k)^2}}{x_0^Tv_1}
\end{align*}
due to $t < m$. We obtain
\begin{equation*}
    \frac{1 - (y_{t+1}^Tx^*)^2}{(y_{t+1}^Tx^*)^2} = \frac{{\sum_{k=2}^d (x_{t+1}^Tv_k)^2}}{(x_{t+1}^Tv_1)^2} 
    \leq 
    \frac{{\sum_{k=2}^d (x_{0}^Tv_k)^2}}{(x_{0}^Tv_1)^2}
    = \frac{1 - (y_0^Tx^*)^2}{(y_0^Tx^*)^2}
\end{equation*}
and we finally have $\Delta_{t+1} = 1-y_{t+1}^Tx^* \leq 1 - y_0^Tx^* = \Delta_0$.
\end{proof}

\begin{proof}[Proof of Lemma~\ref{lemma:recurrence-iteration}]
By Lemma~\ref{lemma:xt-vk-eq}, we have
\begin{align*}
    x_{t+1}^T v_k &= \left(1-\eta+ \eta (\lambda_k + (\lambda^* - \lambda_1) \mathds{1}_{k = 1}) \right) x_t^Tv_k  +  \frac{1}{2} \eta \| x_t \| (y_t - x^*)^T F_k(\hat{y}_t^1,\cdots,\hat{y}_t^d) (y_t - x^*) \\
    & \qquad + \eta \big( G_{S_t}(\bar{y}_t^1,\cdots,\bar{y}_t^d) (x_t-(x_t^Ty_0)y_0) \big)^T v_k.
\end{align*}
Since $S_t$ is sampled uniformly at random, we have $E[f_{S_t}(y)]=f(y)$ for all $y \in \mathbb{R}^d$, which leads to
\begin{align*}
    E[G_{S_t}(\bar{y}_t^1,\cdots,\bar{y}_t^d)] = E[E[G_{S_t}(\bar{y}_t^1,\cdots,\bar{y}_t^d)| \bar{y}_t^1, \bar{y}_t^2, \cdots, \bar{y}_t^d ]] = 0.
\end{align*}
Therefore,
\begin{equation}
\label{proof:x-t-1-v-1-exp}
\begin{aligned}
    & E[(x_{t+1}^Tv_1)^2 \, | \, x_t]
    = \big[ (1-\eta+ \eta \lambda^*) x_t^Tv_1 +  \frac{1}{2} \eta \| x_t \| (y_t - x^*)^T F_1(\hat{y}_t^1,\cdots,\hat{y}_t^d) (y_t - x^*)  \big]^2 \\
    & \qquad \qquad + \eta^2 (x_t-(x_t^Ty_0)y_0)^T E[G_{S_t}(\bar{y}_t^1,\cdots,\bar{y}_t^d)^T v_1 v_1^T G_{S_t}(\bar{y}_t^1,\cdots,\bar{y}_t^d)] ( x_t-(x_t^Ty_0)y_0).
\end{aligned}
\end{equation}
In the same way, for $2 \leq k \leq d$, we have
\begin{equation}
\label{proof:x-t-1-v-k-exp}
\begin{aligned}
    & E[(x_{t+1}^Tv_k)^2 \, | \, x_t]
    = \big[ (1-\eta+ \eta \lambda_k) x_t^Tv_k +  \frac{1}{2} \eta \| x_t \| (y_t - x^*)^T F_k(\hat{y}_t^1,\cdots,\hat{y}_t^d) (y_t - x^*) \big]^2 \\
    & \qquad \qquad + \eta^2 (x_t-(x_t^Ty_0)y_0)^T E[G_{S_t}(\bar{y}_t^1,\cdots,\bar{y}_t^d)^T v_k v_k^T G_{S_t}(\bar{y}_t^1,\cdots,\bar{y}_t^d)] (x_t-(x_t^Ty_0)y_0).
\end{aligned}
\end{equation}
Using the definition of $M$ and $\| \sum_{k=1}^d v_k v_k^T \| = 1$, we have
\begin{equation}
\label{proof:w-t-norm-3}
\begin{aligned}
   \eta^2 ( x_t-(x_t^Ty_0)y_0 )^T
   \sum_{k=1}^d
   E[\| G_{S_t}(\bar{y}_t^1,\cdots,\bar{y}_t^d)^T  v_k \|^2 ] ( x_t-(x_t^Ty_0)y_0 ) \leq \eta^2 K \| x_t-(x_t^Ty_0)y_0 \|^2.
\end{aligned}
\end{equation}
Using \eqref{proof:x-t-1-v-1-exp}, \eqref{proof:x-t-1-v-k-exp}, \eqref{eq:xt-norm-first-bound}, \eqref{eq:xt-vk-second-bound}, \eqref{proof:w-t-norm-3} and the Cauchy-Schwarz inequality for the cross term as
\begin{equation}
\label{proof:cauchy-schwarz-derivation}
\begin{aligned}
    \frac{1}{2} \eta \| x_t \| & \sum_{k=1}^K (1-\eta+\eta (\lambda_k + (\lambda^* - \lambda_1) \mathds{1}_{k=1} )) x_k^Tv_k (y_t - x^*)^T F_1(\hat{y}_t^1,\cdots,\hat{y}_t^d) (y_t - x^*) \\
    & \qquad \qquad \leq \frac{1}{2} \eta M (1-\eta+\eta \lambda^*) \| x_t \| \| y_k - x^* \|^2,
\end{aligned}
\end{equation}
we have
\begin{equation}
\label{proof:w-t-norm-combined}
\begin{aligned}
    E[\|x_{t+1}\|^2 | x_t] & \leq (1-\eta+\eta \lambda^*)^2 \| x_t \|^2 + \frac{1}{2} \eta M (1-\eta+\eta \lambda^*) \| x_t \| \| y_k - x^* \|^2 \\ 
    & \qquad + \frac{1}{4} \eta^2 M^2 \|x_t \|^2 \| y_k - x^* \|^4 + \eta^2 K \| x_t - (x_t^Ty_0) y_0 \|^2.
\end{aligned}
\end{equation}
%Using the Cauchy–Schwarz inequality, we have
%\begin{equation}
%\label{proof:w-t-norm-CS}
%\begin{aligned}
%    \sum_{k=1}^d | x_t^Tv_k| \left | (y_t - x^*)^T F_k(\hat{y}_t^1,\cdots,\hat{y}_t^d) (y_t - x^*) \right |
%    &\leq \sqrt{\sum_{k=1}^d (x_t^Tv_k)^2} \sqrt{\sum_{k=1}^d \big[ (y_t - x^*)^T F_k(\hat{y}_t^1,\cdots,\hat{y}_t^d) (y_t - x^*) \big]^2} \\
%    &\leq {M} \|x_t\| \|y_t - x^* \|^2.
%\end{aligned}
%\end{equation}
Using $\| x_t - (x_t^Ty_0) y_0 \|^2 \leq \|x_t \|^2$ in \eqref{proof:w-t-norm-combined}, we obtain
\begin{equation}
\label{proof:w-t-norm-final}
\begin{aligned}
    E[\|x_{t+1}\|^2 | x_t] & \leq \Big[ \big( 1-\eta+\eta \lambda^* + \frac{1}{2} \eta {M} \| y_t - x^* \|^2 \big)^2 + \eta^2 K \Big] \|x_t\|^2 \\
    &= \Big[ \big( 1-\eta+\eta \lambda^* + \eta {M} (1-y_t^Tx^*) \big)^2 + \eta^2 K \Big] \|x_t\|^2,
\end{aligned}
\end{equation}
which establishes the first statement.

In the same way, using \eqref{proof:x-t-1-v-k-exp}, \eqref{eq:xt-vk-first-bound}, \eqref{eq:xt-vk-second-bound}, \eqref{proof:w-t-norm-3} and the Cauchy-Schwarz inequality similarly to \eqref{proof:cauchy-schwarz-derivation}, we have
\begin{equation}
\label{proof:x-t-1-v-k-exp-sum-2}
\begin{aligned}
    E \Big[ \sum_{k=2}^d (x_{t+1}^T v_k)^2 | x_t \big]
    & \leq
    \bigg[ (1-\eta+\eta \bar{\lambda}) \sqrt{\sum_{k=2}^d (x_t^Tv_k)^2} + \frac{1}{2} \eta {M} \|x_t \| \|y_t - x^*\|^2 \bigg]^2 \\
    &\quad + \eta^2 K \|  x_t-(x_t^Ty_0)y_0 \|^2.
\end{aligned}
\end{equation}
By Lemma~\ref{lemma:y0-St-necessary-condition}, we have $\Delta_t \leq \Delta_0 \leq 1-1/\sqrt{2}$ and thus $y_t^Tx^* \geq 1/\sqrt{2}$ and $y_0^Tx^* \geq 1/\sqrt{2}$. Since $y_t^Tx^* \geq 0$, using \eqref{eq:yt-opt-gap-bound}, we have
\begin{align*}
    \frac{1}{2} \eta {M} \|x_t \| \|y_t - x^*\|^2 
    \leq \eta M  \sqrt{\Delta_t} \sqrt{\sum_{k=2}^d (x_t^Tv_k)^2}.
\end{align*}
As a result of \eqref{eq:yt-y0-bound} which we can use since $\Delta_t \leq \Delta_0$, we obtain
\begin{equation}
\label{proof:x-t-1-v-k-exp-sum-final}
\begin{aligned}
    E \Big[ \sum_{k=2}^d (x_{t+1}^T v_k)^2 | x_t \big] 
    & \leq
    \Big[ 1-\eta+\eta \bar{\lambda} + \eta M \sqrt{\Delta_t} \Big]^2 \sum_{k=2}^d (x_{t}^T v_k)^2 + 8 \eta^2 K \|x_t\|^2  \sum_{k=2}^d (y_0^Tv_k)^2,
\end{aligned}
\end{equation}
which shows the second statement in the lemma.

Lastly, from \eqref{proof:x-t-1-v-1-exp}, we have
\begin{align*}
    E[(x_{t+1}^Tv_1)^2 | x_t] 
    &\geq \left[ (1-\eta+\eta \lambda^*) x_t^Tv_1 + \frac{1}{2} \eta \| x_t \| (y_t-x^*)^T F_1(\hat{y}_t^1,\cdots,\hat{y}_t^d) (y_t-x^*) \right]^2
\end{align*}
By \eqref{eta-condition-1} and \eqref{eta-nu-leq-1}, we have $\eta ( 1 - \lambda^* + M \Delta_0 \sqrt{2}) \leq 1$.
Since $1/(1-\Delta_0) \leq \sqrt{2}$ by \eqref{Delta-0-condition}, we further have
\begin{equation*}
    \eta \left( \frac{M \Delta_0}{1-\Delta_0} + 1 - \lambda^* \right) \leq 1.
\end{equation*}
Due to $\Delta_t \leq \Delta_0$, this implies that
\begin{align*}
    (1-\eta+\eta \lambda^*) x_t^Tv_1 - \frac{1}{2} \eta M \| x_t \| \| y_t-x^* \|^2
    &= \Big[ \left( 1 - \eta + \eta \lambda^* \right) (1-\Delta_t) - \eta M \Delta_t \Big] \|x_t\| \\
    &= \Big[ 1 - \eta \Big( \frac{M \Delta_t}{1-\Delta_t} + 1 - \lambda^* \Big) \Big] (1-\Delta_t) \| x_t \| \\
    &\geq \Big[ 1 - \eta \Big( \frac{M \Delta_0}{1-\Delta_0} + 1 - \lambda^* \Big) \Big] (1-\Delta_t) \| x_t \| \\
    &\geq 0.
\end{align*}
Since $(a+b)^2 \geq (a-c)^2$ holds if $a \geq c$ and $|b| \leq c$, we finally have
\begin{align*}
    E[(x_{t+1}^Tv_1)^2 | x_t]  &\geq \Big[ (1-\eta+\eta \lambda^*) x_t^Tv_1 - \frac{1}{2} \eta M \| x_t \| \| y_t-x^* \|^2 \Big]^2 \\
    &= \left[ 1-\eta+\eta \lambda^* - \eta M \left( \frac{1-y_t^Tx^*}{y_t^Tx^*} \right)  \right]^2 (x_t^Tv_1)^2 \\
    &= \Big[ \alpha(\eta) - \frac{\eta M \Delta_t}{1-\Delta_t} \Big]^2 (x_{t}^T v_1)^2.
\end{align*}
%Lastly, from \eqref{proof:x-t-1-v-1-exp} and \eqref{eq:variance-bound}, we have
%\begin{equation*}
%\begin{aligned}
%    E[(x_{t+1}^Tv_1)^2 \, | \, x_t]
%    &\leq \Big( (1-\eta+ \eta \lambda^*) x_t^Tv_1 +  \frac{1}{2} \eta M \| x_t \| \| y_t - x^* \|^2 \Big)^2 + 4 \eta^2 K \left[ \sum_{k=2}^d (x_t^Tv_k)^2 + \|x_t\|^2 (1-(y_0^Tx^*)^2 \right] \\
%    &\leq \left( 1-\eta+ \eta \lambda^* +  \eta M \left( \frac{1-y_t^Tx^*}{y_t^Tx^*} \right) \right)^2 (x_t^Tv_1)^2 + 4 \eta^2 K \left[ \sum_{k=2}^d (x_t^Tv_k)^2 + \|x_t\|^2 (1-(y_0^Tx^*)^2 \right] \\
%\end{aligned}
%\end{equation*}
\end{proof}

\begin{proof}[Proof of Lemma~\ref{lemma:recurrence-epoch}]
By Lemma~\ref{lemma:y0-St-necessary-condition}, we have $\Delta_t \leq \Delta_0$. Repeatedly applying Lemma~\ref{lemma:recurrence-iteration}, we have
%and using the fact that $\|x_0\|=1$, we have
\begin{equation}
    \label{eq:norm-square-bound}
    \begin{aligned}
    E[ \|x_{t}\|^2 | x_0] 
    = E[ E[ \|x_t\|^2 | x_{t-1} ] | x_0] 
    &\leq \big[ \big( \alpha(\eta) + \eta {M} \Delta_0 \big)^2 + \eta^2 K \big] E[ \|x_{t-1} \|^2 | x_0] \\
    & \leq \big[ \big( \alpha(\eta) + \eta {M} \Delta_0 \big)^2 + \eta^2 K \big]^t \|x_0\|^2.
    \end{aligned}
\end{equation}
Using \eqref{eq:norm-square-bound}, we have
\begin{equation}
\label{eq:xt-vk-square-bound-exp}
\begin{aligned}
    E \Big[ \|x_t\|^2 \sum_{k=2}^d (y_0^Tv_k)^2 \Big] 
    &= E \Big[ E \Big[ \|x_t\|^2 \sum_{k=2}^d (y_0^Tv_k)^2 | x_0 \Big] \Big]
    = E \Big[ E \big[ \|x_t\|^2  | x_0 \big] \sum_{k=2}^d (y_0^Tv_k)^2 \Big] \\
    &= E \Big[ \big[ ( \alpha(\eta) + \eta {M} \Delta_0)^2 + \eta^2 K \big]^t \|x_0\|^2 \sum_{k=2}^d (y_0^Tv_k)^2 \Big] \\
    &= \big[ ( \alpha(\eta) + \eta {M} \Delta_0)^2 + \eta^2 K \big]^t 
    E \Big[  \sum_{k=2}^d (x_0^Tv_k)^2 \Big].
\end{aligned}
\end{equation}
Using Lemma~\ref{lemma:recurrence-iteration} and that $\Delta_t \leq \Delta_0$, we have
\begin{equation}
\label{eq:xt-vk-square-Bt-exp}
    E \Big[ \sum_{k=2}^d (x_{t}^T v_k)^2 \big] \leq
    \big( \beta(\eta) + \eta M \sqrt{\Delta_0}  \big)^2 E \Big[ \sum_{k=2}^d (x_{t-1}^T v_k)^2 \Big] + 8\eta^2 K E \Big[ \|x_{t-1}\|^2  \sum_{k=2}^d (y_0^Tv_k)^2  \Big].
\end{equation}
By induction on \eqref{eq:xt-vk-square-Bt-exp} using \eqref{eq:xt-vk-square-bound-exp}, we have
\begin{equation*}
\begin{aligned}
    E \Big[ \sum_{k=2}^d (x_{t}^T v_k)^2 \big] &\leq
    \big( \beta(\eta) + \eta M \sqrt{\Delta_0}  \big)^2 E \Big[ \sum_{k=2}^d (x_{t-1}^T v_k)^2 \Big] \\
    & \qquad + 8\eta^2 K \big[ ( \alpha(\eta) + \eta {M} \Delta_0 )^2 + \eta^2 K \big]^{t-1} E \Big[ \sum_{k=2}^d (x_0^Tv_k)^2 \Big] \\
    %&\leq E \Big[ \sum_{k=2}^d (x_0^Tv_k)^2 \Big] \Big[ \big( \beta(\eta) + \eta M \sqrt{\Delta_0} \big)^{2t} \\ 
    %& \qquad + 8\eta^2 K \sum_{s=1}^t \big[ \beta(\eta) + \eta M \sqrt{\Delta_0}  \big]^{2(t-s)} \big[ \big( \alpha(\eta) + \eta {M} {\Delta_0} \big)^2 + \eta^2 K \big]^{s-1} \Big] \\
    &\leq E \Big[ \sum_{k=2}^d (x_0^Tv_k)^2 \Big] \Big[ \big( \beta(\eta) + \eta M \sqrt{\Delta_0} \big)^{2t} \\ 
    & \qquad + 8\eta^2 K \sum_{s=1}^t \big( \alpha(\eta) + \eta M \sqrt{\Delta_0}  \big)^{2(t-s)} \big[ ( \alpha(\eta) + \eta {M} \sqrt{\Delta_0})^2 + \eta^2 K \big]^{s-1} \Big] \\
    &\leq E \Big[ \sum_{k=2}^d (x_0^Tv_k)^2 \Big] \Bigg[ \big( \beta(\eta) + \eta M \sqrt{\Delta_0} \big)^{2t} \\
    & \qquad + 8 \big( \alpha(\eta) + \eta {M} \sqrt{\Delta_0} \big)^{2t} \Bigg[ \Bigg( 1 + \frac{\eta^2 K}{( \alpha(\eta) + \eta {M} \sqrt{\Delta_0})^2} \Bigg)^{t} - 1 \Bigg] \Bigg] \\
    &\leq E \Big[ \sum_{k=2}^d (x_0^Tv_k)^2 \Big] \bigg[ \big( \beta(\eta) + \eta M \sqrt{\Delta_0} \big)^{2t} \\
    & \qquad + 8 \big( \alpha(\eta) + \eta {M} \sqrt{\Delta_0} \big)^{2t} \Bigg[ \exp \Bigg( \frac{\eta^2 K t}{( \alpha(\eta) + \eta {M} \sqrt{\Delta_0})^2} \Bigg) - 1 \Bigg].
\end{aligned}
\end{equation*}
By \eqref{eta-condition-2} and \eqref{eta-nu-leq-1}, we have $\eta \big( 1- \lambda^* - M \sqrt{\Delta_0} + \sqrt{Km} \big) \leq 1$, which leads to
\begin{equation*}
    0 \leq \frac{\eta^2 K t}{\big( \alpha(\eta) + \eta {M} \sqrt{\Delta_0} \big)^2} \leq 1.
\end{equation*}
Using $\exp(x) - 1\leq 2x$ for $x \in [0,1]$, we have
\begin{align*}
    E \Big[ \sum_{k=2}^d (x_{t}^T v_k)^2 \big] &\leq E \Big[ \sum_{k=2}^d (x_0^Tv_k)^2 \Big] \left[ \big( \beta(\eta) + \eta M \sqrt{\Delta_0} \big)^{2t} 
    + 16 \eta^2 K t \big( \alpha(\eta) + \eta {M} \sqrt{\Delta_0} \big)^{2(t-1)} \right].
\end{align*}
On the other hand, using $\Delta_t \leq \Delta_0$ and  Lemma~\ref{lemma:recurrence-iteration}, we have
\begin{equation}
\label{eq:xt-v1-square-bound-exp-Bt}
    E[(x_{t}^Tv_1)^2] = E[ E[(x_{t}^Tv_1)^2 | x_{t-1}]]
    \geq \left[ \alpha(\eta) -  \frac{\eta M\Delta_0}{1-\Delta_0} \right]^2 E[(x_{t-1}^Tv_1)^2].
\end{equation}
By induction on \eqref{eq:xt-v1-square-bound-exp-Bt} using $\Delta_t \leq \Delta_0$, we finally have
\begin{align*}
    E[(x_{t}^Tv_1)^2] 
    \geq \left[ \alpha(\eta) -  \frac{\eta M\Delta_0}{1-\Delta_0} \right]^{2t} E[(x_0^Tv_1)^2].
\end{align*}
%In a similar way, using Lemma~\ref{lemma:recurrence-iteration} for $E [ (x_{t}^T v_1)^2 | x_t]$, $\sum_{k=2}^d (x_s^Tv_k)^2 \leq (x_s^Tv_1)^2$ for $s < t$ and \eqref{eq:norm-square-bound}, we have
%\begin{align*}
%    E[ (x_{t}^Tv_1)^2 | B_t ] &\leq \Big[ \alpha(\eta) + \eta M \Big( \frac{\Delta}{1-\Delta} \Big) \Big]^2 E[(x_{t-1}^T v_1)^2 | B_{t-1}] + 4\eta^2K \Big[ E \big[ \sum_{k=2}^d (x_{t-1}^T v_k)^2 | B_{t-1} \big] + \|x_t\|^2 \sum_{k=2}^d (y_0^Tv_k)^2 \Big] \\
%    &\leq \left( \Big[ \alpha(\eta) + \eta M \Big( \frac{\Delta}{1-\Delta} \Big) \Big]^2 + 4 \eta^2 K \right) E[(x_{t-1}^T v_1)^2 | B_{t-1}] \\
%    & \quad + 4\eta^2 K \big[ \big( \alpha(\eta) + \eta {M} \Delta \big)^2 + \eta^2 K \big]^{t-1} \sum_{k=2}^d (y_0^Tv_k)^2 \\
%    &\leq \left( \Big[ \alpha(\eta) + \eta M \Big( \frac{\Delta}{1-\Delta} \Big) \Big]^2 + 4 \eta^2 K \right)^{t} (y_0^Tv_1)^2 \\ 
%    & \quad + 4\eta^2 K \sum_{s=1}^t \left[ \Big[ \alpha(\eta) + \eta M \Big( \frac{\Delta}{1-\Delta} \Big) \Big]^2 + 4 \eta^2 K \right]^{t-s} \big[ \big( \alpha(\eta) + \eta {M} {\Delta} \big)^2 + \eta^2 K \big]^{s-1} \cdot \sum_{k=2}^d (y_0^Tv_k)^2 \\
%    &\leq \left( \Big[ \alpha(\eta) + \eta M \Big( \frac{\Delta}{1-\Delta} \Big) \Big]^2 + 4 \eta^2 K \right)^{t} (y_0^Tv_1)^2 \\
%    & \quad + 4 \eta^2 K t \sum_{k=2}^d (y_0^Tv_k)^2 \left[ \Big[ \alpha(\eta) + \eta M \Big( \frac{\Delta}{1-\Delta} \Big) \Big]^2 + 4 \eta^2 K  \right]^{t-1}.
%\end{align*}
\end{proof}
\begin{proof}[Proof of Lemma~\ref{lemma:expectation-bound-epoch}]
By \eqref{eta-condition-3} and \eqref{def:nu-4}, we have \eqref{eta-condition-2}. Also, \eqref{eta-condition-3}, \eqref{def:nu-5} and the fact that $\sqrt{2 \Delta_0} \leq 1$ which holds from \eqref{Delta-0-condition} imply \eqref{eta-condition-1}. Therefore, by Lemma~\ref{lemma:recurrence-epoch}, we have
\begin{equation}
\label{delta-m-delta-0-bound}
\begin{aligned}
    \delta_m
    \leq
    \left[
    \left( \frac{\beta(\eta) + \eta M \sqrt{\Delta_0}}{\alpha(\eta) - \eta M \Delta_0/(1-\Delta_0)} \right)^{2m}
    + 
    \frac{16 \eta^2 K m \big[ \alpha(\eta) + \eta {M} \sqrt{\Delta_0} \big]^{2(m-1)}}{[ \alpha(\eta) - \eta M \Delta_0 /(1-\Delta_0) ]^{2m}} 
    \right] 
    \delta_0
\end{aligned}
\end{equation}
where 
$$
\delta_t = \frac{E[\sum_{k=2}^d (x_t^Tv_k)^2]}{E[(x_t^Tv_1)^2]}.
$$
By \eqref{Delta0-upper-bound} which follows from \eqref{Delta-0-condition} and the fact that $(1+x)^m \leq \exp(mx)$ for all $x \in \mathbb{R}$, we have
\begin{align*}
   \left( \frac{\beta(\eta) + \eta M \sqrt{\Delta_0}}{\alpha(\eta) - \eta M \Delta_0/(1-\Delta_0)} \right)^{2m}
   &\leq \left( 1 - \frac{\eta(\lambda^*-\bar{\lambda}-2M\sqrt{\Delta_0})}{1-\eta+\eta(\lambda^*- M \sqrt{\Delta_0})} \right)^{2m} \\
   &\leq \exp \left( - \frac{2\eta m(\lambda^*-\bar{\lambda}-2M\sqrt{\Delta_0})}{1-\eta+\eta(\lambda^*- M \sqrt{\Delta_0})} \right).
\end{align*}
Since \eqref{eta-condition-3}, \eqref{def:nu-5} and \eqref{eta-nu-leq-1} imply 
\begin{equation*}
    \frac{2\eta m (\lambda^*-\bar{\lambda}-2M\sqrt{\Delta_0})}{1-\eta+\eta(\lambda^*- M \sqrt{\Delta_0})} \leq 1,    
\end{equation*}
using the fact that $\exp(-x) \leq 1-x/2$ for $0 \leq x \leq 1$, we have
\begin{equation}
    \label{expectation-bound-1}
    \exp \left( - \frac{2\eta m(\lambda^*-\bar{\lambda}-2M\sqrt{\Delta_0})}{1-\eta+\eta(\lambda^*- M \sqrt{\Delta_0})} \right) \leq 1 - \frac{\eta m(\lambda^*-\bar{\lambda}-2M\sqrt{\Delta_0})}{1-\eta+\eta(\lambda^*- M \sqrt{\Delta_0})} = 1- 2\rho.
\end{equation}
On the other hand, by \eqref{Delta0-upper-bound} and the fact that $(1+x)^n \leq \exp(nx)$, we have
\begin{equation}
\label{expectation-bound-2}
\begin{aligned}
    \frac{16 \eta^2 K m \big[ \alpha(\eta) + \eta {M} \sqrt{\Delta_0} \big]^{2(m-1)}}{[ \alpha(\eta) - \eta M \Delta_0 /(1-\Delta_0) ]^{2m}}
    &\leq
    \frac{16 \eta^2 Km}{(\alpha(\eta)+\eta M \sqrt{\Delta_0})^2}
     \left( 1 + \frac{2\eta M \sqrt{\Delta_0}}{\alpha(\eta)-\eta M \sqrt{\Delta_0}} \right)^{2m} \\
    &\leq \frac{16 \eta^2 Km}{( \alpha(\eta)+\eta M \sqrt{\Delta_0})^2} \exp \left( \frac{4 \eta m M \sqrt{\Delta_0}}{\alpha(\eta)-\eta M \sqrt{\Delta_0}} \right).
\end{aligned}
\end{equation}
By \eqref{eta-condition-3}, \eqref{def:nu-4} and \eqref{eta-nu-leq-1}, we have 
\begin{equation*}
    \eta \left(
1- \lambda^* - M \sqrt{\Delta_0} + \frac{64K}{\lambda^*-\bar{\lambda}-2M\sqrt{\Delta_0}} \right) \leq 1,
\end{equation*}
which leads to
\begin{equation}
\label{expectation-bound-3}
    \begin{aligned}
        \frac{\rho}{2} - \frac{16 \eta^2 Km}{(\alpha(\eta)+\eta M \sqrt{\Delta_0})^2} \geq \frac{\eta m (\lambda^* - \bar{\lambda} - 2M\sqrt{\Delta_0}) }{4\big(1-\eta+\eta(\lambda^*+M\sqrt{\Delta_0})\big)}
         - \frac{16 \eta^2 Km}{(1-\eta+\eta(\lambda^* + M \sqrt{\Delta_0}))^2}
         \geq 0.
    \end{aligned}
\end{equation}
In a similar way, by \eqref{eta-condition-3}, \eqref{def:nu-5} and \eqref{eta-nu-leq-1}, we have 
\begin{equation*}
    \eta \left( 1-\lambda^*+M\sqrt{\Delta_0}+ \frac{4mM \sqrt{\Delta_0}}{\log 2} \right) \leq 1,
\end{equation*}
which results in
\begin{equation}
\label{expectation-bound-4}
    \begin{aligned}
        \exp \left( \frac{4 \eta m M \sqrt{\Delta_0}}{\alpha(\eta)-\eta M \sqrt{\Delta_0}} \right) \leq 2.
    \end{aligned}
\end{equation}
Using \eqref{expectation-bound-1}, \eqref{expectation-bound-2}, \eqref{expectation-bound-3} and \eqref{expectation-bound-4} in \eqref{delta-m-delta-0-bound}, we finally have 
\begin{align*}
    \frac{E[ \sum_{k=2}^d (x_m^Tv_k)^2]}{E[(x_m^Tv_1)^2]} 
    &
    \leq
    (1-\rho) \cdot \frac{E[ \sum_{k=2}^d (x_0^Tv_k)^2]}{E[(x_0^Tv_1)^2]}.
\end{align*}
\end{proof}
\begin{proof}[Proof of Theorem~\ref{thm:convergence}]
Since $\eta$, $s$ and $x_0 = \tilde{x}_0$ satisfy
\eqref{Delta-0-condition}, \eqref{L-condition} (or \eqref{eta-condition}) and \eqref{eta-condition-3}, by Lemmas~\ref{lemma:y0-St-necessary-condition} and ~\ref{lemma:expectation-bound-epoch}, we have 
$$
\tilde{\Delta}_1 = \Delta_m \leq \Delta_0 = \tilde{\Delta}_0, \quad \tilde{\delta}_1 = \delta_m \leq (1-\rho) \delta_0 = (1-\rho) \tilde{\delta}_0.
$$
By repeatedly applying the same argument, we have
$\tilde{\delta}_{\tau} \leq (1-\rho)^{\tau} \tilde{\delta}_0$.
Since $\tau \geq ({1}/{\rho}) \log (\tilde{\delta}_0/\epsilon)$, we finally obtain
\begin{equation*}
    \tilde{\delta}_{\tau} \leq (1-\rho)^{\tau} \tilde{\delta}_0
    \leq \exp(-\tau \rho) \tilde{\delta}_0
    \leq \epsilon.
\end{equation*}
This completes the proof.
\end{proof}

\end{document}